\documentclass[12pt]{amsart}
\usepackage{fullpage}
\usepackage[usenames,dvipsnames,svgnames,table]{xcolor}    
\usepackage{pgfplots}
\usepackage{graphicx}
\usepackage{tikz}
\usepackage{tikz-cd}
\usetikzlibrary{decorations.pathreplacing}
\usepackage{lineno, hyperref}
\usepackage{amssymb,amsfonts, amsmath, amsthm}
\usepackage{epstopdf}
\usepackage[shortlabels]{enumitem}
\usepackage{ulem}
\usepackage{cancel}
\usepackage{thmtools} 
\usepackage{mathtools} 
\usepackage{mathrsfs}
\usepackage{yhmath}
\usepackage{longtable}
\usepackage{makecell}
\usepackage{multirow}

\title{Classification of Lipschitz unimodal function germs}

\newtheorem{thm}{Theorem}[section]
\newtheorem{lem}[thm]{Lemma}
\newtheorem{prop}[thm]{Proposition}
\newtheorem{defn}[thm]{Definition}
\newtheorem{cor}[thm]{Corollary}

\newtheorem{question}[thm]{Question}

\numberwithin{equation}{section}

\newcommand{\bb}{\mathbb}
\newcommand{\al}{\mathcal}
\newcommand{\ak}{\mathfrak}

\newcommand{\supp}{{\rm supp}}

\newcommand{\fil}{{\rm fil }}

\newcommand{\Lmod}{{\rm Lmod}}
\newcommand{\smod}{{\rm Smod}}
\newcommand{\jac}{{\rm Jac}}
\newcommand{\Lip}{{\rm Lip}}
\newcommand{\hot}{{\rm h.o.t}}

\author{Nhan Nguyen, Maria Ruas and Saurabh Trivedi}
\address{FPT University, Danang, Vietnam}
\email{nguyenxuanvietnhan@gmail.com}
\address{ICMC, University of Sao Paulo, Sao Carlos, Brazil}
\email{maasruas@icmc.usp.br}
\address{Indian Institute of Technology Goa, at GEC Campus, Farmagudi, Ponda Goa 403110}
\email{saurabh@iitgoa.ac.in}
 

\begin{document}

\setlength{\parindent}{0pt}
\parskip .2cm
\maketitle
\begin{abstract} In this paper, we introduce the notion of Lipschitz modality for isolated singularities 
$ f: (\mathbb{C}^n, 0) \to (\mathbb{C}, 0)$ 
and provide a complete classification of Lipschitz unimodal singularities of corank~2 with non-zero $4$-jets. 
As a consequence, such singularities are Lipschitz unimodal  if they deform to $J_{10}$ but not to $J_{3,0}$. 
Furthermore, we show that singularities with vanishing $6$-jets have Lipschitz modality at least~$2$, thus establishing an upper bound for the order of Lipschitz unimodality.
\end{abstract}



\normalem 

\section{Introduction}

In Singularity Theory, a fundamental achievement in the study of singularities of function germs was the classification of elementary catastrophes and singularities of mappings by Thom and Arnold during the 1960s and 1970s. For isolated singularities, Arnold's work \cite{Arnold74}, \cite{Arnold76}, (see also \cite{Arnold3}) provides a complete classification of simple, unimodal, and bimodal singularities in both real and complex case. Similar classifications have been extended to positive characteristics by Greuel and Nguyen \cite{Greuel-Duc} for simple singularities and by Nguyen \cite{NHD} for unimodal and bimodal singularities.

In Lipschitz classification, two function germs are considered equivalent if there is a bi-Lipschitz homeomorphism in the source transforming one germ into the other. Unlike the topological setting, bi-Lipschitz equivalence admits moduli, as shown by Henry and Parusi\'nski \cite{Henry2, Henry3}. They proved that the Lipschitz type of germs in the family  
$$f_t: (\bb C^2, 0) \to (\bb C, 0), \ f_t(x, y) = x^3 + txy^4 + y^6$$ 
varies continuously in the parameter $t$, which implies that in any small neighborhood of a given $f_t$, there are infinitely many distinct Lipschitz types. Motivated from this result, in a previous work \cite{nsr}, we classified Lipschitz simple germs in the complex case. We defined a germ to be Lipschitz simple if  its $k$-jet  has a neighbourhood (for sufficiently large $k$) intersecting only finitely many Lipschitz classes. An interesting consequence of this classification was the observation that a germ is Lipschitz simple if and only if it does not deform to $J_{10}$ whose nonquaratic part is defined as $ f_t(x,y) = x^3 + txy^4 + y^6$. We also showed in \cite{NS} that this classification agrees with the classification of Lipschitz simple germs under left-right bi-Lipschitz equivalence.

It is worth noting that the bi-Lipschitz equivalence of complex analytic set germs admits no moduli, as proved by Mostowski \cite{Mostowski}. This result was later extended to subanalytic sets by Parusiński \cite{Parusinski}, more generally to polynomially bounded o-minimal structures by Valette and the first author \cite{valette-nhan} (see also \cite{Valette2023}) , and to the p-adic setting by Halupczok and Yin \cite{Halupczok-Yin}.

Since the composition of germs of bi-Lipschitz homeomorphism with analytic function germs is not a group action, we do not have a nice definition of modality as in the smooth case.  In this article we introduce the notion of Lipschitz modality which appears to be natural for bi-Lipschitz equivalence (see Definition \ref{def_Lipschitz_modality}), intuitively, a germ $f: (\mathbb{C}^n, 0) \to (\mathbb{C}, 0)$ has Lipschitz modality at least $m$ if, for sufficiently large $k$, its $k$-jet lies in the closure of a $2m$-dimensional semialgebraic subset of the $k$-jet space whose elements are pairwise non-Lipschitz equivalent. We  give a classification of complex Lipschitz unimodal germs of corank 2 with non-zero 4-jets, extending our earlier work on Lipschitz simple germs. Although the method can, in principle, be applied to all corank 2 germs and potentially to higher corank cases, we restrict our attention to those with non-zero 4-jets in order to take advantage of Arnold’s smooth classification, which is available up to this order. This significantly reduces the amount of computations required for the classification.

One of the crucial steps of our classification is proving that germs in the family  
$$x^3 + b x^2 y^3 + y^9 + cxy^7 + z_1^2 + \dots + z_{n-2}^2$$ have Lipschitz modality $2$ (see Section \ref{section3}). In fact, these are Lipschitz bimodal singularities with the smallest possible Milnor number. Since the uniqueness of the Thom Splitting Lemma is not known in the context of bi-Lipschitz equivalence, this result does not follow directly from the finding of the first author in \cite{nhan2}, which showed that the Lipschitz type of germs in the family
$$ x^3 + b x^2 y^3 + y^9 + cxy^7$$ varies continuously with the parameters $(b, c)$. 
We show this by studying higher-order invariants of the Henry–Parusiński type that depend non-trivially on both parameters 
$b$ and 
$c$ using ideas from \cite{nhan2} (see also the recent work of Migus, Paunescu, and Tibar \cite{Migus}) and arguments similar to those in \cite[Section 3]{nhan1}. 

Based on Arnold's list of classification and the method developed in \cite{nhan1}, we introduce a systematic approach to the classification of Lipschitz unimodal germs with non-zero 4-jets. The classification proceeds through the following steps:

\textbf{Step 1.} Showing that the germs in $J_{3,0}$ whose non-quadratic part is given by 
$$f_{b,c}(x, y) = x^3 + b x^2 y^3 + y^9 + cxy^7$$ 
have Lipschitz modality $2$. This implies that all function germs deforming to $J_{3,0}$ must have Lipschitz modality of at least $2$. This is proved in Section 3.

\textbf{Step 2.} Proving that every germ with a filtration greater than or equal to $9$ (with respect to the weight $(3,1)$), necessarily deforms to $J_{3,0}$. By Step 1, such germs have Lipschitz modality at least 2 and are therefore excluded from our classification. This substantially reduces the set of potential candidates for Lipschitz unimodal singularities. We also prove that the remaining singularities do not deform to $J_{3,0}$. This is done in Section \ref{section4}.

\textbf{Step 3.} We develop a technique for verifying the Lipschitz triviality of a deformation using the Newton diagram. This is presented in Section~\ref{section5}.

\textbf{Step 4.} After Step 2, most singularities are eliminated; only a few families remain to be checked for their Lipschitz modality. From Step (3), we know that these families possess at least $(m-1)$ Lipschitz-trivial directions, where $m$ is the number of parameters in the family. Moreover, germs in these families  do not deform to $J_{3,0}$. Using inductive arguments, we obtain a complete classification of corank-$2$ Lipschitz unimodal germs with nonzero $4$-jets . Furthermore, we prove that all singularities with vanishing $6$-jets deform to $J_{3,0}$, and hence have Lipschitz modality at least $2$. This provides an upper bound for the order of Lipschitz unimodal germs. The detailed results and the final list of Lipschitz unimodal germs are given in Section \ref{section7}.

We conclude the paper with some open questions in Section \ref{section7}.

Throughout the paper, we denote by $\mathbb{C}^* = \mathbb{C} \setminus {0}$. For $x\in \bb R^n$, we denote by $B(x, \varepsilon)$ the open ball radius $\varepsilon$ centered at the $x$. For $X \subset \mathbb{R}^n$, we denote by  $\overline{X}$ the closure of $X$ in $\mathbb{R}^n$. To compare the asymptotic behavior of functions $\varphi$ and $\psi$ near 0, we employ the standard notations $\varphi = o(\psi)$ (or $\varphi \ll \psi$)  and $\varphi = O(\psi)$. 
Given non-negative functions $f, g: X \to \mathbb{R}$, we write $f\lesssim g$ if there exists a positive constant $C$ such that $f(x) \leq C g(x)$ for all $x \in X$; and $f \sim g$ if $f\lesssim g$ and $g\lesssim f$.

\section{Definitions and preliminary results}\label{section2}
We denote by $\mathcal{E}_n$ the $\mathbb{C}$-algebra of complex analytic function germs $f : (\mathbb{C}^n, 0) \to \mathbb{C}$. The \emph{maximal ideal} $\mathfrak{m}_n \subset \mathcal{E}_n$ consists of those germs $f \in \mathcal{E}_n$ such that $f(0) = 0$. 

Let $\mathcal{R}_n$ be the group of germs of biholomophic maps $\varphi: (\mathbb{C}^n, 0) \to (\mathbb{C}^n, 0)$, with the group operation given by composition. The group $\mathcal{R}_n$ acts on $\mathcal{E}_n$ by composition, given by
$$
\varphi \cdot f = f \circ \varphi.
$$
Two germs $f, g \in \mathfrak{m}_n$ are said to be \emph{smoothly right equivalent}, denoted $f \sim_{\mathcal{R}} g$, if they lie in the same orbit under this action; that is, if there exists $\varphi \in \mathcal{R}_n$ such that $f = g \circ \varphi$.

Given $f \in \al E_n$, we denote by $\jac(f)$ the ideal generated by the partial derivatives of $f$, called the Jacobian ideal of $f$. The codimension of a germ $f$ is defined to be the dimension of $\mathbb{C}$-vector space $\al E_n/ \jac(f)$ which is also called the\textit{ Milnor number}, denoted $\mu(f)$. It is known that the Milnor number is a topological invariant in the complex case as proved by Milnor \cite{Milnor}.

We denote by $J^k(n, 1)$ the set of $k$-jets at $0 \in \bb C^n$ of germs in $\al E_n$, it is a $\mathbb{C}$-vector space isomorphic to the vector space of all polynomials in $(x_1, \ldots, x_n)$ with degree $\leq k$. Let $J^k_0(n, 1)$ denote the set of $k$-jets at $0$ of germs in $\ak m_n$. It is obvious that $J^k_0(n, 1)$ is a vector subspace of $J^k(n, 1)$.

Note that a nonsingular germ (i.e., $f \in {\ak m}_n \setminus {\ak m}_n^2$) is equivalent to the projection onto the first coordinate, $x_1$, by the submersion theorem. In what follows, we consider only germs in $\ak m_n^2$, which are often referred to as hypersurface singularities. 

A germ $f\in \ak m_2^n$ is called $k$-determined if for any $g\in \ak m_2^n$ with $j^k f = j^k g$ we have $f\sim_{\al R} g$; $f$ is called finitely determined if it is $k$-determined for some $k \in \bb N$.

In fact, a germ $f$ is finitely determined if and only if it has an isolated singularity at $ 0 $. Moreover, it is known that if $ f $ has an isolated singularity at $ 0 $, then it is $ (\mu(f) + 1) $-determined (see, for example, \cite{Greuel2007}). Now assume that $ f $ is $ k $-determined. Since the Milnor number is upper semicontinuous, every germ in a neighborhood of $ j^k f(0) $ has Milnor number less than or equal to $ \mu(f) $, and hence is also $ (\mu(f) + 1) $-determined. For convenience, from now on we will call any integer $ k \ge \mu(f) + 1 $ \emph{sufficiently large} for $ f $.

For real analytic function germs, being an isolated singularity is not sufficient for finite determinacy, for example $f(x,y) = (x^2 + y^2)^2$ has an isolated singularity in the real case but does not have finite codimension and therefore it is not finitely determined.

The corank of a germ $f\in \ak m_n^2$ is defined as 
$$corank (f) = n - rank (Hess(f))$$ where $Hess(f)$ is the Hessian matrix of $f$ at $0$. Note that singularities of corank $0$ are equivalent to quadratic forms by the Morse lemma. 

For germs of non-zero corank, we have the Thom splitting lemma \cite{Thom75} (see also \cite{Greuel2007},\cite{Pfister2000}), which says:

If a germ $f\in \ak m_n^2$ with an isolated singularity at $0$ has corank $c$, then there exists a germ $g \in \ak m_c^2$ such that 
$$f(x_1, \ldots, x_n) \sim_{\al R} g(x_1, \ldots, x_c) + x_{c+1}^2 + \dots + x_n^2.$$
Moreover, $g$ is uniquely determined up to a diffeomorphism, that is, if $$g_1 (x_1, \ldots, x_c) + x^2_{c+ 1} + \ldots + x^2_n \sim_{\al R} g_2 (x_1, \ldots, x_c) + x^2_{c+ 1} + \ldots + x^2_n $$ then $g_1 \sim_{\al R} g_2$. 
As a result, while writing a germ in $m_n^2$, we often omit the quadratic part and write only its non-quadratic component. For example by writing  $J_{10}: x^3 + txy^4 + y^6$ we mean $J_{10}$ is a class of germs of the forms $x^3 + txy^4 + y^6 + z_1^2 + \ldots + z_{n-2}^2$. 
The uniqueness of Thom splitting lemma remains unknown for bi-Lipschitz equivalence.

For a germ $ f \in \mathfrak{m}_n^2 $, the \emph{smooth modality} of $ f $, denoted by $ \smod(f) $, is defined as the smallest integer $ m $ such that a neighborhood of the $ k $-jet of $ f $ at 0, $ j^k(f)(0) $, for sufficiently large $ k $, can be covered by $ m $ parametric families of orbits under the action of $ \mathcal{R}_n^k $, the group of $ k $-jets of diffeomorphisms from $ (\mathbb{C}^n, 0) $ to itself.

In fact, the smooth modality of $ f $ can be defined using a Rosenlicht stratification under the action of $ \mathcal{R}_n^k $, where $ k $ is taken sufficiently large for $ f $ (see \cite{Greuel-Duc}). This stratification partitions the jet space $ J^k_0(n, 1) $ into locally closed, Zariski-constructible subsets $ \{X_1, \ldots, X_s\} $ such that:
\begin{enumerate}
    \item[(i)] each $ X_i $ is invariant under the action of $ \mathcal{R}_n^k $,
    \item[(ii)] the orbit space $ X_i / \mathcal{R}_n^k $ is an algebraic variety, and
    \item[(iii)] the natural projection $ p_i: X_i \to X_i / \mathcal{R}_n^k $ is a surjective morphism.
\end{enumerate}

Then, the smooth modality of $ f $ is defined as
$$
\smod(f) = \max_i \left\{ \dim p_i(X_i \cap U) \right\},
$$
where $ U $ is a sufficiently small neighborhood of $ j^k(f)(0) $.

Germs with smooth modality 0 are called \emph{smooth simple germs}; those with modality 1 and 2 are called \emph{smooth unimodal} and  \emph{smooth bimodal} germs, and so on.

A \emph{deformation} of a germ $ f \in \mathfrak{m}_n $ is an analytic map germ
$$
F: (\mathbb{C}^n \times \mathbb{C}^m, 0) \to (\mathbb{C}, 0), \quad (x, t) \mapsto F_t(x) = F(x, t),
$$
such that $ F_0(x) = f(x) $.

Given $ f, g \in \mathfrak{m}_n $, we say that $ f $ and $ g $ are \emph{bi-Lipschitz right equivalent}, denoted $f\sim_{Lip} g$ if there exists a germ of a bi-Lipschitz homeomorphism 
$$
\varphi: (\mathbb{C}^n, 0) \to (\mathbb{C}^n, 0)
$$
such that $ f = g \circ \varphi $. Note that if $ \varphi $ is a diffeomorphism, then $ f $ and $ g $ are \emph{smoothly right equivalent}; hence, smooth right equivalence is stronger than bi-Lipschitz right equivalence.

We then define the Lipschitz modality as follows:

\begin{defn}\label{def_Lipschitz_modality} \rm
Let $ f \in \mathfrak{m}_n^2 $ be a finitely determined germ and $k\in \bb N$ be sufficiently large for $f$.  Let $ m $ be the largest integer for which there exists a semialgebraic set  $S \subset J_0^k(n,1)$ of dimension $m$ such that the following conditions hold:
\begin{enumerate}
    \item $ j^kf (0) \in \overline{S} $;
    \item For any distinct $g_1, g_2 \in S$, the germs  at the origin corresponding to $g_1, g_2$ are not bi-Lipschitz right equivalent.
\end{enumerate}
The \emph{Lipschitz modality} of $ f $ is given by
$$
\Lmod (f) = \bigg\lfloor \frac{m}{2} \bigg\rfloor.
$$
\end{defn}

In \cite{nsr}, we defined a finitely determined germ 
$ f \in \mathfrak{m}_n^2 $ to be \emph{Lipschitz simple} if there exist 
only finitely many orbits in a small neighborhood of $ j^k f(0) $ in 
$ J^k_0(n,1) $ for sufficiently large $ k $. It is immediate that if 
$ f $ is Lipschitz simple, then $ \Lmod(f) = 0 $. A natural question 
is whether the converse holds. The answer is yes: from the classification 
in \cite{nsr}, we know that if $ f $ is not Lipschitz simple, then it 
deforms into the class $ J_{10} $, which contains germs with 
$ \Lmod \geq 1 $. Therefore, if $ \Lmod(f) = 0 $, then $ f $ must 
be Lipschitz simple in the sense of \cite{nsr}. 

We adopt the following terminology: $ f $ is \emph{Lipschitz simple} 
if $ \Lmod(f) = 0 $, \emph{Lipschitz unimodal} if $ \Lmod(f) = 1 $, 
\emph{Lipschitz bimodal} if $ \Lmod(f) = 2 $, and so on. Thus, the notion 
of Lipschitz modality introduced in this article generalizes the concept of 
Lipschitz simple germs defined in \cite{nsr}.

Comparing with the notion of smooth modality, one immediately obtains that 
    $$\Lmod(f)\leq \smod(f).$$
Let $ \mathcal{D} $ be a family of function germs. The germ $ f $ is said to \emph{deform to} $ \mathcal{D} $, denoted $ f \to \mathcal{D} $, if there exists a deformation $ F(x, t) $ of $ f $ such that for some sufficiently small $ t \neq 0 $, the germ $ F_t(x) $ is smoothly right equivalent to some germ in $ \mathcal{D} $.

A family $ \mathcal{C} $ is said to \emph{deform to} a family $ \mathcal{D} $, written $ \mathcal{C} \to \mathcal{D} $, if every germ in $ \mathcal{C} $ deforms to $ \mathcal{D} $.

It follows from the definition that if a germ $ f $ deforms to a family $ \mathcal{D} $ of germs of smooth (resp. Lipschitz) modality $ m $, then $ f $ has smooth (resp. Lipschitz) modality at least $ m $.

Arnold \cite{Arnold76} provided a complete classification of germs in 
$ \mathfrak{m}_n^2 $ with isolated singularities and non-zero $4$-jets. 
This classification includes smooth simple, unimodal, and bimodal germs, 
with each class represented by a normal form denoted by 
a letter together with a subscript indicating its Milnor number. For example, 
the simple singularities are
\begin{center}
	\footnotesize
	\begin{tabular}{c|c|c}
		\hline 
		Name & Normal form & $\mu(f)$ \\ \hline
		$A_k$ & $x^{k+1}$ & $k \geq 1$\\
		$D_k$ & $x^2y + y^{k-1}$ & $k\geq 4$\\
		$E_6$ & $x^3 + y^4$ & 6\\
		$E_7$ & $x^3 + xy^3$  & 7\\
		$E_8$ & $x^3 + y^5 $ & 8\\  \hline
	\end{tabular}\\ 
    \vspace{.2cm}
    Smooth Simple Singularities
\end{center}

Let us  list corank $2$ singularities as classified by Arnold \cite{Arnold76} relevant for this paper. Throughout this section, $ a = a_0 + \dots + a_{k-2} y^{k-2}$ for $k>1$ and $ a = 0$ for $k = 1$;

\subsection{The corank 2 singularities with nonzero 3-Jets} These include smooth simple germs $A, D, E_6, E_7, E_8$ and the following singularities:  

\begin{center}
\footnotesize
\begin{longtable}{|c|c|c|c|c|}
\caption{Non-simple corank $2$ singularities with nonzero $3$-jets}
\label{table_corank2}\\
\hline
Name & Smooth normal form & Restrictions &  $\mu(f)$ & $\smod(f)$ \\ 
\hline
$J_{k,0}$        & $x^3 + bx^2 y^k + y^{3k} + {c} xy^{2k+1}$  & $k >1, 4 b^3 + 27 \neq 0$ & $ 6k - 2$ & $k -1$\\
$J_{k,i}$      & $x^3 + x^2 y^k + {a} y^{3k+ i} $ &  $k >1, i >0, a_0 \neq 0$ & $6k + i- 2$ & $k - 1$\\
$E_{6k}$        & $x^3 + y^{3k + 1} + {a}x y^{2k+1}$  & $k \geq 1$ & $6k$& $k-1$ \\
$E_{6k+1}$        & $x^3 + xy^{2k + 1} + { a} y^{3k+2}$  & $k \geq 1$ &   $6k + 1$&  $k-1$\\
$E_{6k+2}$        & $x^3 +  y^{3k + 2} + { a} x y^{2k+2}$  & $k \geq 1$ & $6k +2$& $k-1$ \\ \hline
\end{longtable}
\end{center}
Here, 
${c} =c_0 + \dots + c_{k-3}y^{k-3}$ for $k>2$ and ${c} = 0$ for $k = 2$.

\subsection{The corank 2 singularities with zero $3$-jets and nonzero $4$-jets}
These include singularities of classes $X$, $Y$, $Z$ and $W$ which are described as follows:

\textbf{Classes $X$ and $Y$}:
\begin{center}
\footnotesize
\begin{longtable}{|c|c|c|c|c|}
\caption{ Class $X$ and $Y$ with $k>1$}
\label{table_corank2}\\
\hline
Name & Smooth normal form & Restrictions &  $\mu(f)$ & $\smod(f)$  \\ 
\hline 
$X_{k,0}$ &$x^4+bx^3y^k+ax^2y^{2k} +xy^{3k}$       & $\Delta \neq 0$, $a_0b_0 \neq 9$ & $12k-3$& $3k-2$ \\
 $X_{k,p}$     & $x^4 + ax^3y^k + x^2y^{2k} + by^{4k+p}$& $a_0^2 \neq 4$, $b_0\neq 0$, $p>0$ & $12k-3+p$ & $3k-2$\\
$Y_{r,s}^k$     & $[(x+ay^k)^2+by^{2k+s}](x^2+y^{2k+r})$  & $1\leq s\leq r$,$k>1$,$a_0\neq 0 \neq b_0$  &$12k-3+r+s$ &  $3k-2$\\
  \hline
\end{longtable}
\end{center}
Here, $\Delta = 4(a_0^3+b_0^3)-a_0^2b_0^2 - 18a_0b_0+27$ and $b = b_0+\cdots+b_{2k-2}
y^{2k-2}$.

\begin{center}
\footnotesize
\begin{longtable}{|c|c|c|c|c|}
\caption{ Class $X$ and $Y$ with $k=1$}
\label{table_corank2}\\
\hline
Name & Smooth normal form & Restrictions &  $\mu(f)$ & $\smod(f)$  \\ 
\hline 
$X_{1,0} = X_9$ &$x^4+a_0 x^2 y^2 + y^4$       & $a_0^2 \neq 4$ & $9$& $1$ \\
 $X_{1,p}$     & $x^4 + x^2y^2 + a_0 y^{4+p}$ & $a_0 \neq 0$ & $9+p$ & $1$\\
$Y_{r,s}^1$     & $x^{4+r} + a_0 x^2 y^2 + y^{4+s}$ & $a_0\neq 0$  &$9+r+s$ &  $1$\\
\hline
\end{longtable}
\end{center}

\textbf{Class $Z$:}

For singularities $Z_{i, 0}^k$ and $Z_{\mu}^k$ $(k>1)$ normal forms are $f = (x + ay^k)f_2$ where $a_0\neq 0$ and $f_2$ is given in the following table: 
\begin{center}
\footnotesize
\begin{longtable}{|c|c|c|c|c|}
\caption{Class $Z$ with $k>1$}
\label{table_corank2}\\
\hline
Name & Normal form & Restrictions &  $\mu(f)$ & $\smod(f)$  \\ 
\hline 
$Z_{i,0}^k$ & $x^3+dx^2y^{k+1} + cxy^{2k+2i+1}+ y^{3k+3i}$       & $4d^3+27\neq 0$, $k>1$ $i\geq 0$ & $12k+6i-3$& $3k+i-2$\\
$Z^k_{12k+6i-1}$ & $x^3+bxy^{2k+2i+1}+y^{3k+3i+1}$ & $k>1$,$i\geq0$ & $12k+6i-1$  & $3k+i-2$\\
$Z^k_{12k+6i}$    & $x^3+xy^{2k+2i+1}+by^{3k+3i+2}$  & $k>1$,$i\geq0$ & $12k+6i$& $3k+i-2$\\
$Z^k_{12k+6i+1}$ & $x^3 + bxy^{2k+2i+2} + y^{3k+3i+2}$ & $k>1$,$i\geq0$ & $12k+6i+1$ &$3k+i-2$\\
\hline
\end{longtable}
\end{center}
For $k = 1$, class $Z$ consists of the following families: 
\begin{itemize}
    \item $Z_{i,0}$, $Z_{6i+11}$, $Z_{6i+12}$, $Z_{6i + 13}$ ($i>0$) which have the normal forms $f = yf_2$ where $f_2$ is given in Table \eqref{table_corank2}. 
    \item $Z_{i,p}: y(x^3 + x^2y^{i+1} +  b y^{3i+p+3})$, $b_0\neq 0$, $i >0$, $p>0$
\end{itemize}
where $ b = b_0 + \dots + b_{2k+i-2}y^{2k+i-2}$ and $c = c_0 + \dots + c_{2k+i-3}y^{2k + i -3}$.

\textbf{Class W:}
\begin{center}
\footnotesize
\begin{longtable}{|c|c|c|c|c|}
\caption{Class $W$}
\label{table_corank2}\\
\hline
Name & Smooth normal form & Restrictions &  $\mu(f)$ & $\smod(f)$ \\ 
\hline 
$W_{12k}$ & $x^4 + y^{4k+1} + a x y^{3k+1} + c x^2 y^{2k+1}$       &  & $12k$& $3k-2$\\
$W_{12k+1}$ & $x^4 + xy^{3k+1} + a x^2 y^{2k+1} + c y^{4k+2}$ &  & $12k+1$& $3k-2$\\
$W_{k,0}$    & $x^4 + b x^2 y^{2k+1} + a xy^{3k+2} + y^{4k+2}$  &  $b_0^2 \neq 4$ & $12k+3$& $3k-1$\\
$W_{k,i}$    & $x^4 + a x^3 y^{k+1} +  x^2y^{2k+1} +b y^{4k+2+i}$  & $i>0, b_0 \neq 0$ & $12k+3+i$& $3k-1$\\
$W_{k,2q-1}^{\#}$ & $(x^2 + y^{2k+1})^2 + b x y^{3k+1+q} + a y^{4k+2+q}$ & $q>0, b_0\neq 0$ & $12k+2 + 2q$ &$3k-1$\\
$W_{k,2q}^{\#}$ & $(x^2 + y^{2k+1})^2 + b x^2 y^{2k+1+q} + a x y^{3k+2+q}$ & $q>0, b_0\neq 0$ & $12k+3 + 2q$ &$3k-1$\\
$W_{12k+5}$ & $x^4 + xy^{3k+2} + a x^2 y^{2k+2} + b y^{4k+3}$ &  & $12k+5$ &$3k-1$\\
$W_{12k+6}$ &  $x^4 + y^{4k+3} + a xy^{3k+3} + b x^2y^{2k+2}$ &  & $12k+6$ & $3k-1$\\
\hline
\end{longtable}
\end{center}
Here, $k \geq 1$ and $b = b_0 + \dots + b_{2k-1} y^{2k-1}$, $c = c_0 + \dots + c_{2k-2} y^{2k-2}$.

\subsection{The corank 2 singularities of smooth modality $2$} There are $8$ infinite series and $8$ exceptional families where $a = a_0 + a_1 y$. 
\begin{center}
\footnotesize
\begin{longtable}{|c|c|c|c|}
\caption{Bimodal germs of corank 2: the $8$ infinite series}
\label{table_corank2}\\
\hline
Name & Smooth normal form & Restrictions &  $\mu(f)$ \\ 
\hline 
$J_{3,0}$        & $x^3 +  bx^2y^3 + y^9 +  cx y^{7}$  & $4 b^3 + 27\neq 0$ & $16$  \\
$J_{3,p}$        & $x^3 +  x^2y^3 +  ay^{9+p}$  & $p>0, a_0\neq 0$ & $16+p$  \\
$Z_{1,0}$        & $y(x^3 +  dx^2y^2  +  cxy^5 + y^6)$  & $4 d^3 + 27\neq 0$ & $15$  \\
$Z_{1,p}$        & $y(x^3 +  x^2y^2  +  ay^{6+p})$  & $ p>0, a_0 \neq 0$ & $15+p$ \\
$W_{1,0}$        & $x^4 +  {a}x^2y^3  +  y^{6}$  & $ a_0^2 \neq 4$ & $15$  \\
$W_{1,p}$        & $x^4 + x^2y^3  +  ay^{6+p}$  &$p>0$, $ a_0 \neq 0$ & $15+p$ \\
$W^{\#}_{1,2q-1}$        & $(x^2 +y^3)^2  +  { a}xy^{4+q}$  & $q>0, a_0 \neq 0$ & $15+2q-1$  \\
$W^{\#}_{1,2q}$        & $(x^2 +y^3)^2  +  {a}x^2y^{3+q}$  & $q>0, a_0 \neq 0$ & $15+2q$ \\
\hline
\end{longtable}
\end{center}
\vspace{-.5cm}
\begin{center}
\footnotesize
\begin{longtable}{|c|c|c|c|}
\caption{Bimodal germs of corank 2:  the $8$ exceptional families}
\label{table_corank2}\\
\hline
Name & Smooth normal form & Name &  Smooth normal form \\ 
\hline
$E_{18}$        & $x^3 + y^{10} + axy^7$  & $E_{19}$ & $x^3 + xy^7 + a y^{11}$  \\
$E_{20}$        & $x^3 + y^8 + axy^8$  & $Z_{17}$ & $x^3y + y^{11} + axy^6$  \\
$Z_{18}$        & $x^3y +  xy^6 + ay^9 $  & $Z_{19}$ & $x^{3}y + y^9 + axy^7$  \\
$W_{17}$        & $x^4 +  xy^5 + ay^7 $  & $W_{18}$& $x^4 + y^7 + ax^2 y^4$ \\
\hline
\end{longtable}
\end{center}
\vspace{-0.75cm}
\subsection{The corank $2$ Lipschitz simple singularities}

In \cite{nsr} we show that:
\begin{thm}\label{thm_Lipschitz_simple}
Let $f$ be an isolated singularity of corank $\leq 2$. Then, $f$ is Lipschitz simple if and only if it is smoothly equivalent to one of the following germs:
\begin{center}
	\footnotesize
	\begin{longtable}{|c|c|c|c|c|}
    \caption{Corank 1 and Corank 2 Lipschitz simple germs}\\
    \hline
	Name & Smooth normal form & Restrictions &  $\mu(f)$ & Corank \\ \hline
		$A_k$        & $x^{k+1}$     &  $k\geq 1$             & $k$ & $1$\\\hline 
		$D_k$        & $x^2y +y^{k-1}$  &         $k\geq 4$  & $k$ & \\
	$E_6$        & $x^3 + y^4$  &               &  $6$& \\
		$E_7$        & $x^3 + xy^3$     &                 & $7$& \\
		$E_8$        & $x^3 + y^5 $  &        & $8$& $2$\\
		$X_9$        & $x^4+y^4 + tx^2y^2$              & $t\neq 0$              & $9$& \\
		$T_{2,4,5}$    & $x^4 + y^5 + tx^2y^2$        & $t\neq 0$        & $10$ & \\
		$T_{2,5,5}$	 & $ x^5 +y^5 +  tx^2 y^2 $	      & $    t\neq 0 $	       & $11$ & \\
		$Z_{11}$	 & $x^3 y + y^5 +  txy^4$	   &  &  $11$ &\\	
		$W_{12}$     & $x^4  + y^5 + tx^2y^3$       &         & $12$ & \\
	\hline
	\end{longtable}
\end{center}
\end{thm}

Moreover, we have:
\begin{thm}\label{thm_Lipschitz_simple2}
Let $f$ be an isolated singularity. Then, $f$ is  Lipschitz modal if and only if $f$ deforms to $J_{10}: x^3 + txy^4 + y^6$.
\end{thm}

\section{Lipschitz modality of  $J_{3,0}$}\label{section3}

In this section we will prove that the germs in $J_{3,0}:  x^3 + bx^2y^3 + y^9 + cxy^7, \frac{4}{27}b^3 + 1 \neq 0$ has Lipschitz modality $2$.

Let us start with some notation. Let $f, g: (\bb C^n, 0) \to (\bb C, 0)$ be analytic function germs. Suppose that they are bi-Lipschitz equivalent. Then, there is a bi-Lipschitz homeomorphism $\varphi: (\bb C^n, 0) \to (\bb C^n, 0)$ such that $f = g\circ \varphi$. Let $L$ be the bi-Lipschitz constant of $\varphi$. It has been known (see for example \cite[Lemma 4.5]{nsr}) that 
\begin{equation}\label{equ_J30_1}
L^{-1}\|\nabla g(\varphi(x))\| \leq \|\nabla f(x))\| \leq L \|\nabla g(\varphi(x))\|.
\end{equation}

Given $\delta> 1$, we define
$$V^\delta (f) =\{x \in \bb C^n : f(x)\neq 0, \delta^{-1} \|x\|\|\nabla f(x)\| \leq \|f(x)\|\leq  \delta \|x\|\|\nabla f(x)\|\}.$$
It follows from (\ref{equ_J30_1}) that
\begin{equation}\label{equ_J_30_2}
V^{L^{-2} \delta} (g) \subset \varphi(V^\delta(f)) \subset V^{L^{2}\delta}(g).
\end{equation}

Given $\sigma >0$, $d >0$, we define 
$$ W^{\sigma}(f, d) = \{x \in \bb C^n: |f(x)| \leq \sigma \|x\|^d\}.$$
This implies
\begin{equation}\label{equ_J_30_3}
W^{\sigma L^{-d}}(g, d) \subset \varphi(W^\sigma (f,d)) \subset W^{\sigma L^d} (g, d).
\end{equation}
We define $$\Omega(f, \delta, \sigma,d) = V^{\delta}(f) \cap W^{\sigma}(f, d).$$ 
It follows from (\ref{equ_J_30_2}) and (\ref{equ_J_30_3}) that
\begin{equation}\label{equ_J30_4}
 \Omega(g, \delta L^{-2}, \sigma L^{-d},d) \subset \varphi(\Omega(f,\delta,  \sigma,d)) \subset \Omega(g, \delta L^{2}, \sigma L^{d},d). 
\end{equation}

Now we consider germs in $J_{3,0}$ which are of the following form
$$f_{b,c}(x, y, z) = x^3 + bx^2y^3 + y^9 + cxy^7 + z_1^2+ \dots + z_{n-2}^2$$
 where $b \neq \frac{4}{27}b^3+ 1$, $z = (z_1, \ldots z_{n-2}) \in \bb C^{n-2}$. For convenience, we put $z^2  = z_1^2 + \dots + z^2_{n-2}$. The family $f_{b,c}$ then can be  shorten as
 $$f_{b,c}(x,y, z) = x^3 + bx^2y^3 + y^9 + cxy^7 +z^2.$$
 
 Fix $(b, c) \in \bb C^2$ and fix constants $\delta>1$, $\sigma>0$. Note that $\Omega(f_{b,c}, \delta,\sigma, d)$ is a semialgebraic subset of $\bb C^n \equiv \bb R^{2n}$. Consider the case $d = 9$. We have the following properties:

\begin{lem}\label{lem_3.1}  The tangent cone at $0$ of $\Omega(f_{b,c}, \delta,\sigma,9)$ is contained in the $y$-axis.
\end{lem} 

\begin{proof} Let $v = (m_1, m_2, m_3) \in \bb C^{2}\times \bb C^{n-2}$ be a unit vector of the tangent cone of $\Omega(f_{b,c}, \delta,\sigma,9)$ at $0$. It suffices to show that $m_1 = 0$ and $m_3 = 0$. 

By the curve selection lemma, there exists a real analytic arc $\gamma:[0,\varepsilon) \to \Omega(f_{b,c}, \delta,\sigma,9)$ such that $\gamma(0)=0$ and $v = \lim_{t\to 0}\frac{\gamma(t)}{\|\gamma(t)\|}$. Write $\gamma(t) = (\gamma_x(t), \gamma_y(t), \gamma_z(t))$. By reparameterizing $\gamma$, we may assume that $\|\gamma(t)\| \sim |t|$.  

 Recall that $$\nabla f_{b,c}(x, y, z)= ( 3x^2 + 2bxy^3 + cy^7, 3bx^2y^2 + 9 y^8 + 7 cxy^6, 2z)$$

If $m_3 \neq 0$, we have $|\gamma_x(t)| \sim  |t|^\alpha$, $|\gamma_y(t)| \sim  |t|^\beta$ and $|\gamma_z(t)| \sim  |t|$ where $\alpha, \beta \geq 1$ (put $|t|^\infty = 0$ as a convention). This implies $|f_{b,c}(\gamma(t))| \sim  |t|^2$. Hence $\gamma(t) \not\subset W^{\sigma}(f_{b,c}, 9)$.

If $m_3 = 0, m_1\neq 0$, we have $|\gamma_x(t)| \sim  |t|$,  $|\gamma_y(t)|\sim |t|^\alpha$  and $|\gamma_z(t)|\sim |t|^\beta$ where $\alpha \geq 1$ and $\beta > 1$.  The proof is split into two cases.

\textbf{Case 1:} $2\beta \leq 3$ 

Then $ f_{b,c}(\gamma(t)) \lesssim |t|^{2\beta}$ and $|\nabla f_{b,c}(\gamma(t))| \sim |t|^\beta$. 
Thus, 
$$\|\gamma(t)\|\|\nabla f_{b,c}(\gamma(t))\| \sim |t|^{\beta +1} \gg  |t|^{2\beta} \gtrsim | f_{b,c}(\gamma(t))| \text{ (since } \beta >1)$$ This implies $\gamma$ is not contained in $V^{\delta}(f_{b,c})$.\\

\textbf{Case 2:}  $2\beta > 3$ 

Then, $|f_{b,c}(\gamma(t))| \sim |t|^{3}$, hence yielding that  $$\gamma(t) \not \subset W^{\sigma}(f_{b,c}, 9).$$

All the above cases give contradictions, so $m_1 =  $ and $m_3$ must be $0$.
\end{proof}

\begin{lem}\label{lem_3.2}
On the germ at $0$ of $\Omega(f_{b,c}, \delta, \sigma, 9)$,  we have
\begin{itemize}
    \item either (i) $x = -\frac{2b}{3} y^3 + \frac{c}{2b}y^4 + o(y^4)$ or (ii)  $x = -\frac{c}{2b} y^4 + o(y^4)$.
     \item $z = O(y^8)$.
 
\end{itemize}
\end{lem}

\begin{proof}
We have 
$$\partial f_{b,c} /\partial x =  3x^2 + 2bxy^3 + cy^7,$$ 
$$\partial f_{b,c} /\partial y =  3bx^2y^2 + 9 y^8 + 7 cxy^6$$  and $$\partial f_{b,c} /\partial z = 2z.$$ 
If $(x, y, z) \in \Omega(f_{b,c}, \delta,\sigma,9)$, then  
$$\|\nabla f_{b,c} (x,y, z)\| \|(x, y, z)\| \sim | f_{b,c}(x, y, z)|$$ and
 $$|f_{b,c}(x, y, z)| \lesssim \|(x, y, z)\|^9.$$ 
 
 By Lemma \ref{lem_3.1}, we have $\|(x, y, z)\| \sim |y|$. This implies that $$\|\nabla f_{b,c}(x, y, z)\| \lesssim  \|(x,y,z)\|^8 \sim |y|^8.$$
 Hence, 
$$
\begin{cases}
|3x^2 + 2bxy^3 + cy^7| \lesssim |y|^8,\\
|3bx^2y^2 + 9 y^8 + 7 cxy^6| \lesssim |y|^8 \\
2|z|  \lesssim |y|^8
\end{cases}
$$
It is clear that the conclusion follows from the first and the third inequalities.
\end{proof}

\begin{thm}\label{prop_j30_modal}  All germs in $J_{3, 0}$ have Lipschitz modality $2$.
\end{thm}

\begin{proof} Recall that germs in $J_{3,0}$ are of the following form: 
$$f_{b,c} (x, y, z) = x^3 + bx^2y^3 + cxy^7 + z^2 \frac{4}{27}b^3 + 1 \neq 0,$$
where $z  = (z_1, \ldots, z_{n-2})\in \bb C^{n-2}$ and $z^2 = \sum_i z_i^2$.

Given $(b, c)$ with $b, c\neq 0$, $\frac{4}{27}b^3 + 1 \neq 0$, and constants $\delta > 1$ and $\sigma > 0$, consider the germ at $0$  of $\Omega (f_{b,c}, \delta, \sigma, 9)$.
By Lemma \ref{lem_3.2}, this germs can be separated into two germs, denoted by,  $\Omega_1(f_{b,c}, \delta, \sigma, 9)$ and $\Omega_2(f_{b,c}, \delta, \sigma, 9)$, where points in $\Omega_1(f_{b,c}, \delta, \sigma, 9)$ are of the form 
$$
x = -\frac{2b}{3} y^3 + \frac{c}{2b}y^4 + o(y^4), \quad z = O(y^8),
$$
and points in $\Omega_2(f_{b,c}, \delta, \sigma, 9)$ are of the form 
$$
x = -\frac{c}{2b} y^4 + o(y^4), \quad z = O(y^8).
$$
Note that $\Omega_1(f_{b,c}, \delta, \sigma, 9) \setminus \{0\}$ and $\Omega_2(f_{b,c}, \delta, \sigma, 9) \setminus \{0\}$ are disjoint. In addition, the restrictions of $f_{b,c}$ to $\Omega_1(f_{b,c}, \delta, \sigma, 9)$ and $\Omega_2(f_{b,c}, \delta, \sigma, 9)$, respectively, are 
\begin{equation}\label{equJ_30_3.1}
    f_{b,c}|_{\Omega_1} (x,y,z) = \left( \frac{4}{27} b^3 + 1\right) y^9 - \frac{2bc}{3}y^{10} +  o(y^{10})
\end{equation}
\begin{equation}\label{equJ_30_3.2}
  f_{b,c}|_{\Omega_2} (x, y, z) = y^9   + o(y^{10}).
\end{equation}
Note that the expressions in \eqref{equJ_30_3.1} and \eqref{equJ_30_3.2} are independent of $(\delta, \sigma)$.

Now consider germs $f_{b_1, c_1}$ and $f_{b_2, c_2}$ with $(b_1, c_1)\neq (b_2, c_2)$, $b_i \neq 0$, $ \frac{4}{27}b_i^3+ 1\neq 0$, $c_i\neq 0$ for $i = 1, 2$. Suppose that $f_{b_1, c_1}$ and $f_{b_2, c_2}$ are bi-Lipschitz right equivalent. Let $h: (\bb C^{n}, 0) \to (\bb C^{n},0)$  be a bi-Lipschitz homeomorphism such that $f_{b_1, c_1} =  f_{b_2, c_2} \circ h$ and let $L\geq 1$ be the bi-Lipschitz constant of $h$. By (\ref{equ_J30_4}), we have 
\begin{equation}\label{equJ_30_3.3}
 \Omega( f_{b_2, c_2}, \delta L^{-2}, \sigma L^{-9},9)  \subset h( \Omega(f_{b_1, c_1}, \delta, \sigma,9))\subset \Omega(f_{b_2, c_2},\delta L^{2}, \sigma L^{9},9).
\end{equation}

Let
$$ u_k = \left(\frac{-2b_1}{3}y_k^3 + \frac{c_1}{2b_1}y_k^4 + o(y_k^4), y_k, 0\right)$$
and 
$$\tilde{u}_k = \left(-\frac{c_1}{2b_1}y_k^4 + o(y_k^4), y_k, 0\right)$$
be sequences contained in $\Omega_1(f_{b_1, c_1}, \delta, \sigma, 9)$ and $\Omega_2(f_{b_1, c_1}, \delta, \sigma, 9)$, respectively, tending to $0$ as $k \to \infty$.

Set $v_k = h(u_k)$ and $\tilde{v}_k = h(\tilde{u}_k)$. By \eqref{equJ_30_3.3},  $(v_k)$ and $(\tilde{v}_k)$ are both contained in $\Omega(f_{b_2, c_2}, \delta L^2, \sigma L^9, 9)$ and also tend to $0$ as $k \to \infty$.

The proof is now divided into several cases.

\underline{\textbf{Case 1:}} $(v_k)$ and $(\tilde{v}_k)$ are both in $\Omega_{1}(f_{b_2, c_2}, \delta L^2, \sigma L^9, 9)$

We can write 
$$ v_k = \left(\frac{-2b_2}{3}Y_k^3 + \frac{c_2}{2b_2}Y_k^4 + o(Y_k^4), Y_k, 0\right)$$
and 
$$\tilde{v}_k = \left(\frac{-2b_2}{3}{\tilde{Y}}_k^3+ \frac{c_2}{2b_2}{\tilde{Y}}_k^4 + o({\tilde{Y}}_k^4), {\tilde{Y}}_k, 0\right)$$
where $Y_k$ and $\tilde{Y}_k$ tend to $0$ as $k$ tends to $\infty$.
It follows from \eqref{equJ_30_3.1} that 
\begin{equation}\label{equJ_30_3.4}
    f_{b_2,c_2}(v_k) = \left( \frac{4}{27} b_2^3 + 1\right) Y_k^9 - \frac{2b_2c_2}{3}Y_k^{10} + o(Y_k^{10})
\end{equation}
and 
\begin{equation}\label{equJ_30_3.5}
    f_{b_2,c_2}(\tilde{v}_k) = \left( \frac{4}{27} b_2^3 + 1\right) {\tilde{Y}}_k^9 - \frac{2b_2c_2}{3}{\tilde{Y}}_k^{10}  + o({\tilde{Y}}_k^{10})
\end{equation}
By definition,  $|u_k| \sim |\tilde{u}_k|\sim |y_k|$. Since $h$ is bi-Lipschitz, $|v_k| \sim |\tilde{v}_k|\sim  |u_k|$.  Since $|v_k| \sim |Y_k|$ and $|\tilde{v}_k| \sim |\tilde{Y}_k|$,  so $|y_k| \sim |Y_k| \sim |{\tilde{Y}}_k|$. 

Since $| u_k - \tilde{u}_k| \sim |y_k|^3$, we have $|v_k - \tilde{v}_k| \sim |y_k|^3$. Then,
\begin{equation}
    |Y_k - \tilde{Y}_k| \lesssim  |v_k - \tilde{v}_k|  \sim  |y_k|^3 \sim  |Y_k|^3
\end{equation} 
Hence, 
\begin{equation}\label{equJ_30_3.6}
    \tilde{Y}_k = Y_k + O(Y_k^3)
\end{equation}  
In \eqref{equJ_30_3.5} replacing $\tilde{Y}_k $ with $Y_k + O(Y_k^3)$ we get 
\begin{equation} \label{equJ_30_3.7}
    f_{b_2,c_2}(\tilde{v}_k) = \left( \frac{4}{27} b_2^3 + 1\right) {Y}_k^9 - \frac{2b_2c_2}{3}{Y}_k^{10} + o({Y}_k^{10})
\end{equation}
Recall that
\begin{equation}\label{equJ_30_3.8}
    f_{b_1,c_1} (u_k) = \left( \frac{4}{27} b_1^3 + 1\right) y_k^9 - \frac{2b_1c_1}{3}y_k^{10} + o(y_k^{10})
\end{equation}
\begin{equation}\label{equJ_30_3.9}
  f_{b_1,c_1}(\tilde{u}_k)  = y_k^9 + o(y_k^{10})
\end{equation}
Since $ f_{b_1, c_1} (u_k) = f_{b_2, c_2} (v_k)$ and $f_{b_1, c_1} (\tilde{u}_k) = f_{b_2, c_2} (\tilde{v}_k)$, we have 
\begin{equation}\label{equJ_30_3.10}
    f_{b_1, c_1} (u_k) - f_{b_1, c_1} (\tilde{u}_k)  = f_{b_2, c_2} (v_k) - f_{b_2, c_2} (\tilde{v}_k) 
\end{equation}
Equivalently, $$ \frac{4}{27}b_1^3 y_k^9 -\frac{2b_1c_1}{3}y_k^{10} + o(y_k^{10}) = o(Y_k^{10})$$
Since $|y_k| \sim |Y_k|$, we have $ b_1 = 0$. This contradicts our assumption that $b_i \neq 0$ for $i = 1, 2$.

\underline{\textbf{Case 2:}} $(v_k)$ and $(\tilde{v}_k)$ are both in $\Omega_{2}(f_{b_2, c_2}, \delta L^2, \sigma L^9, 9)$

We write
$$ v_k = \left(-\frac{c_2}{2b_2}Y_k^4 + o(Y_k^4), Y_k, 0\right)$$
and 
$$ \tilde{v}_k = \left(-\frac{c_2}{2b_2}\tilde{Y}_k^4 + o(Y_k^4), \tilde{Y}_k, 0\right)$$
where $Y_k$ and $\tilde{Y}_k$ tend to $0$ as $k$ tends to $\infty$.
By \eqref{equJ_30_3.2},
\begin{equation}
    f_{b_2,c_2}(v_k) = Y_k^9 + o(Y_k^{10})
\end{equation}
and 
\begin{equation}\label{equJ_30_3.11}
    f_{b_2,c_2}(\tilde{v}_k) =  \tilde{Y}_k^9 + o(\tilde{Y}_k^{10})
\end{equation}
By \eqref{equJ_30_3.6}, we can write
$$f_{b_2,c_2}(\tilde{v}_k) = Y_k^9  + o(Y_k^{10})$$
Then, \eqref{equJ_30_3.10} is equivalent to $$\frac{4}{27}b_1^3 y_k^9 -\frac{2b_1c_1}{3}y_k^{10} + o(y_k^{10})  = o(Y_k^{10}).$$
Since $|y_k| \sim |Y_k|$, this implies  that $b_1=0$, and again contradicts the assumption that $b_i \neq 0$ for $i = 1, 2$.

\underline{\textbf{Case 3:}} $(v_k)  \subset \Omega_{1}(f_{b_2, c_2}, \delta L^2, \sigma L^9, 9)$ and $ (\tilde{v}_k) \subset  \Omega_{2}(f_{b_2, c_2}, \delta L^2, \sigma L^9, 9)$. 

We can write 
$$ v_k = \left(\frac{-2b_2}{3}Y_k^3 + \frac{c_2}{2b_2}Y_k^4 + o(Y_k^4), Y_k, 0\right)$$
and 
$$ \tilde{v}_k = \left(-\frac{c_2}{2b_2}\tilde{Y}_k^4 + o(Y_k^4), \tilde{Y}_k, 0\right)$$
where $Y_k$ and $\tilde{Y}_k$ tend to $0$ as $k$ tends to $\infty$.
Thus, 
 \begin{equation}
    f_{b_2,c_2}(v_k) = \left( \frac{4}{27} b_2^3 + 1\right) {Y}_k^9 - \frac{2b_2c_2}{3}{Y}_k^{10} + o({Y}_k^{10})
\end{equation}
and 
\begin{equation}
    f_{b_2,c_2}(\tilde{v}_k) = {\tilde{Y}_k}^9 + o({\tilde{Y}_k}^{10})
\end{equation}
Again, by \eqref{equJ_30_3.6}, we can write
$$ f_{b_2,c_2}(\tilde{v}_k) = {Y_k}^9 + o({Y_k}^{10}).$$
Then, \eqref{equJ_30_3.10} is equivalent to 
\begin{equation}\label{equJ_30_3.12}
    \left( \frac{4}{27} b_1^3 + 1\right) y_k^9 - \frac{2b_1c_1}{3}{y_k}^{10}  + o({y_k}^{10}) =  \left( \frac{4}{27} b_2^3 + 1\right) {Y_k}^9 - \frac{2b_2c_2}{3}{Y_k}^{10} + o({Y_k}^{10})
\end{equation}
and 
\begin{equation}\label{equJ_30_3.13}
    {y_k}^9   + o({y_k}^{10}) = {Y_k}^9   + o(Y_k^{10}).
\end{equation}

Since $|y_k| \sim |Y_k|$, it follows from \eqref{equJ_30_3.12} that 
\begin{equation*}
    \lim_{k\to \infty} \left(\frac{Y_k}{y_k}\right)^9 = \frac{(4/27) b_1^3+1}{(4/27) b_2^3+1}.
\end{equation*}
From \eqref{equJ_30_3.13} we get
\begin{equation}\label{equJ_30_3.14}
    \lim_{k\to \infty} \left(\frac{Y_k}{y_k}\right)^9 = 1.
\end{equation}
It follows that
\begin{equation}\label{equJ_30_3.15}
    b_1^3 = b_2^3.
\end{equation}
On the other hand, by dividing $\frac{4}{27}b_1^3 + 1$ in both sides of \eqref{equJ_30_3.12}, we obtain the following:
\begin{equation}\label{equJ_30_3.16}
  {y_k}^9 -\frac{\frac{2b_1c_1}{3}} {\left( \frac{4}{27} b_1^3 + 1\right)}{y_k}^{10}  + o(y_k^{10}) =   Y_k^9 - \frac{\frac{2b_2c_2}{3}}{   \left( \frac{4}{27} b_1^3 + 1\right)}Y_k^{10} + o(Y_k^{10})
\end{equation}
Subtracting side-by-side \eqref{equJ_30_3.16} by \eqref{equJ_30_3.13}, we have:

\begin{equation}
    \frac{\frac{2b_1c_1}{3}} {\left( \frac{4}{27} b_1^3 + 1\right)} y_k^{10}  = \frac{\frac{2b_2c_2}{3}} {\left( \frac{4}{27} b_1^3 + 1\right)} Y_k^{10} + o(y_k^{10})
\end{equation}
This yields 
\begin{equation}
    \lim_{k\to \infty} \left(\frac{Y_k}{y_k}\right)^{10} = \frac{b_1c_1}{b_2c_2}
\end{equation}
Hence, 
\begin{equation}
    \lim_{k\to \infty} \left(\frac{Y_k}{y_k}\right)^{90} = \left(\frac{b_1c_1}{b_2c_2}\right)^9
\end{equation}
By\eqref{equJ_30_3.14}, we get 
\begin{equation}
    \left(\frac{b_1c_1}{b_2c_2}\right)^9 =1
\end{equation}
Note from \eqref{equJ_30_3.15} that $b_1^3 = b_2^3$. It follows that 
\begin{equation}
    \left(\frac{c_1}{c_2}\right)^9  = 1
\end{equation}
Equivalently,  $$c_1^9 = c_2^9.$$ 
Therefore, in this case, we have $b_1^3 = b_2^3$ and $c_1^9 = c_2^9$.

\underline{\textbf{Case 4:}} $(v_k)  \subset \Omega_{2}(f_{b_2, c_2}, \delta L^2, \sigma L^9, 9)$ and $ (\tilde{v}_k) \subset  \Omega_{1}(f_{b_2, c_2}, \delta L^2, \sigma L^9, 9)$.

Using similar arguments as in Case 3, just interchange the role of $v_k$ and $\tilde{v}_k$, we get 

\begin{equation}
   f_{b_2,c_2} (v_k) =  {Y_k}^9 + o({Y_k}^{10})
\end{equation}
and 
\begin{equation}
   f_{b_2,c_2} (\tilde{v}_k) =  \left( \frac{4}{27} b_2^3 + 1\right) {Y}_k^9 - \frac{2b_2c_2}{3}{Y}_k^{10} + o({Y}_k^{10})
\end{equation}

Then, \eqref{equJ_30_3.10} is equivalent to 
\begin{equation}\label{equJ_30_3.17}
    \left( \frac{4}{27} b_1^3 + 1\right) y_k^9 - \frac{2b_1c_1}{3}{y_k}^{10}  + o({y_k}^{10}) =  {Y_k}^9   + o(Y_k^{10}).
\end{equation}
and 
\begin{equation}\label{equJ_30_3.18}
    {y_k}^9   + o({y_k}^{10}) =  \left( \frac{4}{27} b_2^3 + 1\right) {Y_k}^9 - \frac{2b_2c_2}{3}{Y_k}^{10} + o({Y_k}^{10})
\end{equation}

Since $|y_k| \sim |Y_k|$, it follows from \eqref{equJ_30_3.17} that 
\begin{equation}
    \lim_{k\to \infty} \left(\frac{Y_k}{y_k}\right)^9 = (4/27) b_1^3+1
\end{equation}
From \eqref{equJ_30_3.18} we get
\begin{equation}\label{equJ_30_3.19}
    \lim_{k\to \infty} \left(\frac{Y_k}{y_k}\right)^9 = \frac{1}{(4/27) b_2^3+1}.
\end{equation}
It follows that
\begin{equation}\label{equJ_30_3.20}
    \left(\frac{4}{27} b_1^3+1\right)\left(\frac{4}{27} b_2^3+1\right) = 1
\end{equation}
On the other hand, dividing $\frac{4}{27}b_1^3 + 1$ in both sides of \eqref{equJ_30_3.17} yields
\begin{equation}\label{equJ_30_3.21}
  {y_k}^9 -\frac{\frac{2b_1c_1}{3}} {\left( \frac{4}{27} b_1^3 + 1\right)}{y_k}^{10}  + o(y_k^{10}) =  \frac{1}{   \left( \frac{4}{27} b_1^3 + 1\right)} Y_k^9  + o(Y_k^{10})
\end{equation}
Subtracting side-by-side \eqref{equJ_30_3.21} by \eqref{equJ_30_3.18} gives

\begin{equation}
    \frac{\frac{2b_1c_1}{3}} {\left( \frac{4}{27} b_1^3 + 1\right)} y_k^{10}  = -\frac{2b_2c_2}{3} Y_k^{10} + o(y_k^{10})
\end{equation}
This implies
\begin{equation}
    \lim_{k\to \infty} \left(\frac{Y_k}{y_k}\right)^{10} = -\frac{1}{\frac{4}{27} b_1^3 + 1} \frac{b_1c_1}{b_2c_2}
\end{equation}
Together with \eqref{equJ_30_3.19}, we get 
\begin{equation}
    \left(\frac{b_1c_1}{b_2c_2}\right)^9 = -\left(\frac{4}{27} b_1^3 + 1\right)^{19}
\end{equation}
In summary, for this case we obtain
$$\left(\frac{4}{27} b_1^3+1\right)\left(\frac{4}{27} b_2^3+1\right) = 1$$
and 
$$\left(\frac{b_1c_1}{b_2c_2}\right)^9 = -\left(\frac{4}{27} b_1^3 + 1\right)^{19}.$$

We are now ready to show that the germs in $ J_{3,0} $ have Lipschitz modality $ 2 $. Since these germs have \emph{smooth} modality $ 2 $, they must have \emph{Lipschitz} modality $\leq  2 $. Therefore, it suffices to show that their Lipschitz modality is $\geq 2 $.

Consider the set
$$
V = \left\{(b,c) \in \mathbb{C}^2 \ \middle|\ b \ne 0,\, c \ne 0,\, \frac{4}{27}b^3 + 1 \ne 0,\, \text{and } \left|\frac{4}{27}b^3 + 1\right| \ne 1 \right\}.
$$
It is clear that $ V $ is an open and dense semialgebraic subset of $ \mathbb{C}^2 $. Let $ \mathcal{D} = \{f_{b,c} \mid (b,c) \in V\} $. Every germ in $ J_{3,0} $ deforms to some germ in $ \mathcal{D} $, so it is enough to show that each germ in $ \mathcal{D} $ has Lipschitz modality  $\geq 2 $.

Fix a point $ a_0 = (b_0, c_0) \in V $. We will show that there exists $ \varepsilon > 0 $ such that for any two distinct points $ (b_1, c_1) \ne (b_2, c_2) \in B(a_0, \varepsilon) $, we have $ f_{(b_1, c_1)} \not\sim_{\mathrm{Lip}} f_{(b_2, c_2)} $. This will imply, by definition, that $ f_{b_0, c_0} $ has Lipschitz modality  $ \geq 2 $.

Since $ V $ is open,  we may assume that  $ B(a_0, \varepsilon) \subset V $ when choosing $\varepsilon$ sufficiently small. Suppose, for contradiction, that there exist distinct points $ (b_1, c_1), (b_2, c_2) \in B(a_0, \varepsilon) $ such that $ f_{(b_1, c_1)} \sim_{\mathrm{Lip}} f_{(b_2, c_2)} $. Because $ b_1, b_2 \ne 0 $, neither Case 1 nor Case 2 from above applies. From Case 3 and Case 4, it follows that either 

\begin{equation} \label{equJ_30_3.22}
\begin{cases}
b_1^3 = b_2^3,\\
c_1^9 = c_2^9,
\end{cases}
\end{equation}

or

\begin{equation} \label{equJ_30_3.23}
\begin{cases}
\left(\dfrac{4}{27} b_1^3+1\right)\left(\dfrac{4}{27} b_2^3+1\right) = 1,\\
\left(\dfrac{b_1c_1}{b_2c_2}\right)^9 = -\left(\dfrac{4}{27} b_1^3 + 1\right)^{19}.
\end{cases}
\end{equation}

Now observe that
$$
b_1^3 - b_2^3 = (b_1 - b_2)(b_1^2 + b_1b_2 + b_2^2), \quad
c_1^9 - c_2^9 = (c_1 - c_2)\sum_{i+j=8} c_1^i c_2^j.
$$
We may take $ \varepsilon $ small enough so that both $ |b_1^2 + b_1b_2 + b_2^2| > 0 $ and $ \left| \sum_{i+j=8} c_1^i c_2^j \right| > 0 $; this is possible since $ a_0 \ne (0,0) $. If \eqref{equJ_30_3.22} occurs, then $ b_1 = b_2 $ and $ c_1 = c_2 $, contradicting our assumption that $ (b_1, c_1) \ne (b_2, c_2) $.

Thus, we may assume that only \eqref{equJ_30_3.23} holds.   Observe that $ \left|\dfrac{4}{27} b_0^3 + 1\right| \ne 1 $. Without loss of generality, we may assume that $ \left|\dfrac{4}{27} b_0^3 + 1\right| < 1 $. Then, for sufficiently small $ \varepsilon $, we also have
$$
\left|\dfrac{4}{27} b_1^3 + 1\right| < 1 \quad \text{and} \quad \left|\dfrac{4}{27} b_2^3 + 1\right| < 1.
$$
This implies that the first equation in \eqref{equJ_30_3.23} cannot hold. Therefore, this case also leads to a contradiction.

Hence, $ f_{(b_1, c_1)} \not\sim_{\mathrm{Lip}} f_{(b_2, c_2)} $, completing the proof.
\end{proof}

\section{Germs that deform to $J_{3,0}$}\label{section4}
In this section we classify all corank $2$ germs of non-zero $4$-jets that deform to $J_{3,0}$. 

Since modality is upper semicontinuous and  $J_{3,0}$ contains germs of smooth modality $2$,  any germ deforming to $J_{3,0}$ must have smooth modality at least $2$. According to Arnorld's classificaton, non-zero 4-jets corank 2 germs with smooth modality $\geq 2$ are:

\begin{itemize}
    \item $J_{k,0}$, $J_{k,i}$ with $k \geq 3$, $i > 1$;
    \item $E_{6k}$, $E_{6k+1}$, $E_{6k+2}$ with $k \geq 3$;
    \item Classes $X$ and $Y$ with $k > 1$;
    \item Class $Z$ with $k > 1$; $Z_{i,0}$, $Z_{i,p}$, $Z_{6i+11}$, $Z_{6i+12}$, $Z_{6i+13}$ with $i > 0$, $p > 0$;
    \item $W_{12k}$, $W_{12k+1}$ with $k > 1$; $W_{k,0}$, $W_{k,i}$, $W_{k,2q-1}^{\#}$, $W_{k,2q}^{\#}$, $W_{12k+5}$, $W_{12k+6}$ with $k > 0$.
\end{itemize}

We will prove the following result:

\begin{thm}\label{thm_deform_J30}
    Let $f\in m_n^2$ be a corank $2$ germ  of non-zero $4$-jet with smooth modality $\geq 2$. Then, $f$ does not deform to $J_{3,0}$ if and only if $f$ is smoothly equivalent to one of the germs in the following families: 
     \begin{enumerate}
    \item $Z_{1,0}: x^3 y + sx^2y^3 + txy^6 + y^7$, $4s^3 + 27 \neq 0$
    \item $Z_{1, 1}: x^3 y + x^2 y^3 + sy^8 + t y^9$, $s \neq 0$ 
     \item $Z_{17}: x^3y + y^8 + sxy^6 + txy^7$, $s\neq 0$
    
   \item  $W_{1,0}: x^4 +  sx^2y^3 + t x^2 y^4  +  y^{6}$, $s^2\neq 4$ 
   \item $W_{1,1}: x^4 + x^2y^3  +  sy^{7} + t y^8$, $s \neq 0$
   \item $W_{1,2}: x^4 + x^2y^3  +  sy^{8} + t y^9$, $s \neq 0$
   \item $W_{17}: x^4 + xy^5 + sy^7 + ty^8$, $s\neq 0$
    \item $W_{18}: x^4 + y^7 + sx^2y^4 + tx^2 y^5$, $s\neq 0$
    
    \item $W_{1, 2q-1}^{\#}: (x^2 + y^3)^2 + s x y^{4+q} + t xy^{5 + q}$, $s \neq 0$, $q>0$
    \item $W_{1, 2q}^{\#}: (x^2 + y^3)^2 + sx^2 y^{3+q} + t x^2 y^{4 + q}$, $s \neq 0$, $q >0$
    \end{enumerate}
\end{thm}

Let $w = (w_1,\ldots,w_n) \in \bb Q_{\geq0}^n$ be a non-zero vector and $d \in \bb Q_+$. A polynomial $f(x) = \sum_\alpha c_\alpha x^\alpha$ is called quasihomogenous of type $(w_1,\ldots,w_n; d)$ if for all powers $\alpha$ in $f$:
$$\sum_{i=1}^n w_i \alpha_i = d.$$ Given a monomial $x^{\alpha} = x_1^{\alpha_1}\ldots x_n^{\alpha_n}$, the filtration of $x^{\alpha}$ with respect to weight $w$, is defined as $\fil_w(x^{\alpha}) = \sum_{i=1}^n  w_i \alpha_i$. Then, the filtration of a germ of analytic function $f: (\bb C^n, 0)\to (\bb C, 0)$ is the minimum of the filtrations of the monomials appearing in the Taylor expansion of $f$.

We need the following results to prove Theorem \ref{thm_deform_J30}.

\begin{lem}\label{lem_main_deform}
Let $ w = (3, 1) $ be a weight, and let $ f(x,y) $ be a germ of a polynomial  with an isolated singularity at $ 0 $. If $ \fil_w(f) \geq 9 $, then $ f $ deforms to $ J_{3,0} $.
\end{lem}

\begin{proof}
Define $ g_b(x, y) = x^3 + bx^2 y^3 + y^9 $, where $ 4b^3 + 27 \neq 0 $. Consider the following deformation with generic $ b \in \mathbb{C} $:
$$
F_t(x, y) = f(x, y) + t g_b(x, y).
$$
Since $ \fil_w(f) \geq 9 $, the germ $ F_t $ can be written in the general form:
$$
F_t(x, y) = a_1(t)\, x^3 + a_2(t)\, x^2 y^3 + a_3(t)\, x y^6 + a_4(t)\, y^9 + h_t(x, y),
$$
where $ h_t $ is a family of analytic germs with $ \fil_w(h_t) \geq 10 $ for every $ t $.

For generic $ t $, a suitable change of coordinates of the form $ x \mapsto \alpha_1 x + \alpha_2 y^3 $ and $ y \mapsto \beta y $, one can eliminate the term $xy^6$, that yields
$$
F_t(x, y) \sim_R g_c(x, y) + \widetilde{h}_t(x, y),
$$
for some $ c \in \mathbb{C} $ and some germ $ \widetilde{h}_t $ with $ \fil_w(\widetilde{h}_t) \geq 10 $.

Note that $ g_c $ is quasi-homogeneous of type $ (w; 9) $. A direct computation shows that, for $ c \neq 0 $,
$$
\mathbb{C}\{x, y\}/\jac(g_c) = \langle 1, x, y, xy, y^2, xy^2, y^3, y^4, xy^3, xy^4, y^5, y^6, y^7, y^8, y^9, y^{10} \rangle.
$$
Among these, only $ y^{10} $ has filtration with respect to $ w $ strictly greater than 9. By \cite[Theorem 7.2]{Arnold74}, any germ of the form $ g_c + g' $ with $ \fil_w(g') > 9 $ is smoothly equivalent to $ g_c + a y^{10} $ for some $ a \in \mathbb{C} $. Consequently,
\begin{equation}\label{eq6.1}
F_t \sim_{\al R} g_c + \widetilde{h}_t \sim_{\al R} g_c + a(t)\, y^{10}.
\end{equation}
The germ on the right-hand side of \eqref{eq6.1} belongs to the class $ J_{3,0} $. Therefore, $ f $ deforms to $ J_{3,0} $.
\end{proof}

\begin{thm} \label{thm_deformation}
Germs in families in Theorem \ref{thm_deform_J30} do not deform to $J_{3,0}$.
\end{thm}

\begin{proof}

(i) $Z_{1,0}  \not \to J_{3,0}$ and $W_{1,0}  \not \to J_{3,0}$. 

Let $f$ be a germ in $Z_{1,0}$ or in $W_{1,0}$. Then $\mu(f) = 15$. Since  Milnor number is upper semicontinuous, in a sufficently small neighborhood of $J^k(f)$ for $k$ large enough, there are only germs of Milnor nubmer at most $15$. Since all germs in $J_{3,0}$ have Milnor number $16$, $f$ cannot deform to $J_{3,0}$.

(ii) $Z_{1,1}\not \to J_{3,0}$ (similar arguments can be used to show that $W_{1,1}\not \to J_{3,0}$).

Let $f$ be a germ in $Z_{1,1}$. Then $\mu(f) = 16$. Suppose on the contrary that $f \to J_{3,0}$. Since the germs in $J_{3,0}$ have the same Milnor number $ \mu = 16$, there exists a $\mu $-constant family of function germs $ G_t(x,y) $ such that $ G_0 = f $ and $ G_t \in J_{3,0} $ for $ t $ near $ 0 $.  By \cite[Theorem 1.1]{Bobadilla} ,  the multiplicity $m(G_t)$ must remain constant throughout the family. However, $ m(f) = 4 $, while $ m(G_t) = 3 $ for all $ t \neq 0 $ near $ 0 $ which is a  contradiction. Thus, $Z_{1,1} \not\to J_{3,0}$.

(iii) $ W_{1,2} \not\to J_{3,0} $ (similar arguments can be used to show that $Z_{17},  W_{18}, W^{\#}_{1, 2q-1}, W^{\#}_{1, 2q}\not\to J_{3,0}$). Since $W_{18} \to W_{17}$ and $W_{18} \not \to J_{3,0}$, hence $W_{17} \not \to J_{3,0}.$

Suppose, on the contrary, that there exists a deformation $ \phi_a(x, y) $ such that:
\begin{itemize}
    \item $ \phi_a $ for $a\neq 0$ is of type $ J_{3,0} $:
    $$
    \phi_a(x, y) \sim_{\al R}  x^3 + s x^2 y^3 + y^9 + t x y^7.
    $$
    \item $ \phi_0 $ is of type $ W_{1,2} $:
    $$
    \phi_0(x, y) \sim_{\al R} x^4 + x^2 y^3 + s y^8 + t y^9, \quad s \neq 0.
    $$
\end{itemize}

After a suitable linear change of coordinates, we can assume that the 8-jet of $ \phi_a(x, y) $ has the form:
$$
\phi_a(x, y) = a x^3 + \sum_{i=4}^8 G_i(x, y),
$$
where:
\begin{align*}
G_4(x, y) &= p_0 x^4 + p_1 x^3 y + p_2 x^2 y^2 + p_3 x y^3 + p_4 y^4, \\
G_5(x, y) &= q_0 x^5 + \dots + q_3 x^2 y^3 + q_4 x y^4 + q_5 y^5, \\
G_6(x, y) &= u_0 x^6 + \dots + u_5 x y^5 + u_6 y^6, \\
G_7(x, y) &= v_0 x^7 + \dots + v_7 y^7, \\
G_8(x, y) &= m_0 x^8 + \dots + m_8 y^8.
\end{align*}

Here, $ p_i, q_i, u_i, v_i, m_i $ are smooth functions in $ a $.

From Arnold's classification, for $a\neq 0$ close to $0$:
\begin{itemize}
    \item If $ p_4(a) \neq 0 $, then $ \phi_a $  is of type $ E_6 $.
    \item If $ p_4(a) = 0 $, $ p_3(a) \neq 0 $, then $ \phi_a $ is of type $ E_7 $.
    \item If $ p_4 (a)= p_3(a) = 0 $, $ q_5(a) \neq 0 $, then $ \phi_a $ is of type $ E_8 $.
\end{itemize}
Since $ \phi_a $ is of type $ J_{3,0} $, none of these cases are possible. Therefore, we may assume:
$$
p_4 = p_3 = q_5 \equiv 0.
$$
in a neighborhood of $0$.

To determine the type of $ \phi_0 $, we analyze the coefficients $ p_i, q_i, u_i, v_i, m_i $ at $a = 0$.

\begin{itemize}
    \item If $ p_1(0) \neq 0 $ or $ p_2(0) \neq 0$, the tangent cone of $\phi_0 $ consists of more than one irreducible component, so it cannot be of type $W_{1,2}$.
    \item If $ p_1(0) = p_2(0) = 0 $ and $ p_0(0) = 0  $, the 4-jet of $ \phi_0 $ is zero, hence, $\phi_0$ is not of type $W_{1,2}$.
    \item If $ p_1(0) = p_2(0) = 0 $, $ p_0(0)\neq 0 $, and $ q_4(0)\neq 0$, then $ \phi_0 $ is of type $ W_{13} $.
    \item If $ p_1(0)= p_2(0) =  q_4(0) =  0 $, $ p_0 (0) \neq 0 $, and $ u_6(0)\neq 0 $, then $ \phi_0 $ is of type $ W_{1,0} $, $ W_{1,2q-1}^\# $, or $ W_{1,2q}^\# $.
    \item If $ p_1(0) = p_2(0)=q_4(0)=u_6(0) = 0 $, $ p_0(0) \neq 0 $, and $ u_5(0)\neq 0  $, then $ \phi_0 $ is of type $ Z_{17} $.
    \item If $ p_1(0) = p_2(0)=q_4(0)= u_6(0)=u_5(0)$, $ p_0(0)\neq 0 $, $ q_3 \neq 0 $, and $ v_7 \neq 0 $, then $ \phi_0 $ is of type $ W_{1,1} $.
    \item If  $ p_1(0) = p_2(0)=q_4(0)= u_6(0)=u_5(0) = q_3(0) = 0$, $ p_0(0)\neq  0 $ and $ v_7 (0) \neq 0 $, then $ \phi_0 $ is of type $ W_{18} $.
    \item If  $ p_1(0) = p_2(0)=q_4(0)= u_6(0)=u_5(0) = v_7(0) = 0$, $ p_0(0) \neq 0, q_3(0)\neq 0 $, and $ m_8(0)=0 $, then $ \phi_0 $ is of type $ W_{1,p} $ for $ p > 2 $.
\end{itemize}

Since $ \phi_0 $ is of type $ W_{1,2} $, it follows that 
\begin{equation}\label{equ_deform_0}
\begin{cases}
p_1(0) = p_2(0)=q_4(0)= u_6(0)=u_5(0) = v_7(0) = 0,\\ 
p_0(0) \neq 0, q_3(0) \neq 0 , \ m_8 (0)\neq  0.\\
\end{cases}
\end{equation}

Now, consider $ \phi_a $ with $ a \neq 0 $. By a change of coordinates:
$$
x \mapsto x - \frac{p_2}{3a} y^2, \quad y \mapsto y,
$$
we eliminate the $ x^2 y^2 $ term from $ G_4(x, y) $.

After substitution, to ensure $ \phi_a $ is of type $ J_{3,0} $, the coefficients of $ x y^4 $, $ y^6 $, $ x y^5 $, and $ y^8 $ must vanish. This gives:
\begin{itemize}
    \item The coefficient of $ x y^4 $:
    $$
    q_4 - \frac{p_2^2}{3a} = 0.
    $$
    \item The coefficient of $ y^6 $:
    $$
    u_6 + \frac{2 p_3^2}{27 a^2} - \frac{p_2^3}{9 a^2} = 0.
    $$
    \item The coefficient of $ y^7 $:
    \begin{equation}\label{equ_deform1}
    v_7 - \frac{u_5 p_2}{3a} + \frac{q_3 p_2^2}{9 a^2} - \frac{p_1 p_2^3}{27 a^3} = 0.
    \end{equation}
    \item The coefficient of $ x y^5 $:
    \begin{equation}\label{equ_deform2}
    -\frac{2 q_3 p_2}{3a} + \frac{p_1 p_2^2}{3 a^2} + u_5 = 0.
    \end{equation}
\end{itemize}
We focus on the latter two equations.

\textbf{Claim:} $ \lim_{a \to 0} \frac{p_2}{a} = 0 $.

From \eqref{equ_deform2}, we have:
$$
\frac{p_2}{3a} \left( 2 q_3 - \frac{p_1 p_2}{3a} \right) = u_5.
$$

Since $ \lim_{a\to 0} u_5 = u_5(0) = 0$, we have either
\begin{equation*}
\lim_{a \to 0} \frac{p_2}{a} = 0 \tag{1}
\end{equation*}
or
\begin{equation*}
\lim_{a \to 0} \left( 2 q_3 - \frac{p_1 p_2}{3a} \right) = 0 \tag{2}.
\end{equation*}

If (1) holds, then the claim follows trivially.

If (2) holds, then
\begin{equation}\label{equ_deform3}
\lim_{a \to 0} \frac{p_1 p_2}{3a} = 2 \lim_{a \to 0} q_3.
\end{equation}

From \eqref{equ_deform1},
$$
\frac{p_2}{3a} \left( u_5 - \frac{q_3 p_2}{3a} + \frac{p_1 p_2^2}{9 a^2} \right) = v_7.
$$
As $ a \to 0 $, $ u_5 $ and $ v_7 \to 0 $. Thus, either:
$$
\lim_{a \to 0} \frac{p_2}{a} = 0,
$$
or
\begin{equation}\label{equ_deform4}
\lim_{a \to 0} \left( -\frac{q_3 p_2}{3a} + \frac{p_1 p_2^2}{9 a^2} \right) = 0.
\end{equation}

Again, if the first case holds, then the claim is trivial. We may assume that the second case holds. Substituting \eqref{equ_deform3} into \eqref{equ_deform4}, we can replace $ \frac{p_1 p_2}{3a} $ with $ 2 q_3 $, leading to:
$$
\lim_{a\to 0}\frac{q_3 p_2}{3a} = 0.
$$

Since $ \lim_{a\to 0} q_3 = q_3(0) \neq 0 $, it follows that:
$$
\lim_{a \to 0} \frac{p_2}{a} = 0.
$$
The claim is proved.

Now consider the coefficient of $ y^8 $:
$$
A = m_8 + \frac{u_4 p_2^2}{9 a^2} - \frac{q_2 p_2^3}{27 a^3} + \frac{p_0 p_2^4}{81 a^4} - \frac{v_6 p_2}{3a}.
$$

Since $\lim_{a\to 0} \frac{p_2}{a} = 0 $, all terms involving $ p_2 $ vanish, leaving:
$$
\lim_{a \to 0} A = \lim_{a \to 0} m_8(a) = m_8(0)\neq 0 \text{ (by \eqref{equ_deform_0})}  
$$

The nonzero value $ A $ implies that $ \phi_a $ is of type $ E_{14} $, not $ J_{3,0} $, which is a contradiction. Therefore, $W_{1,2} \not\to J_{3,0}$.
\end{proof}

\begin{proof}[Proof of Theorem \ref{thm_deform_J30}]
By Theorem  \ref{thm_deformation}, germs listed in the statement do not deform to $J_{3,0}$.

It is straightforward to check that that corank-$2$ germs with nonzero $4$-jets and smooth modality $\geq 2$ in Arnold's classification not listed in the statement of the theorem have filtrations $\geq 9$ with respect to the weight $w = (3,1)$. Therefore, by Lemma \ref{lem_main_deform}, these germs deform to $J_{3,0}$. 
\end{proof}

\section{Lipschitz triviality of families}\label{section5}

In this section, we provide a sufficient condition for a family of function germs to be Lipschitz trivial (Theorem~\ref{thm_checking_triviality_1}). This is inspired by the computations used in proving Lipschitz triviality of families in Lemmas~7.12 and~7.14 of~\cite{nsr}. We then apply this result to give a list of Lipschitz trivial families in Theorem \ref{thm_checking_triviality_main}.

\subsection{Thom-Levine's criterion for Lipschitz triviality}

Consider an analytic family of function germs  $$F(x, t) = f_t(x): (\mathbb{C}^n, 0) \to (\mathbb{C}, 0)$$  where $t$ is in a connected open subset $U\subset \bb C$. The family $ F $ is \textit{Lipschitz trivial} over $ U $ if, for each $ t_0 \in U $, there exists a neighborhood $ U_{t_0} \subset U $ of $ t_0 $ and a continuous family of bi-Lipschitz homeomorphisms $ h_t: (\mathbb{C}^n, 0) \to (\mathbb{C}^n, 0) $, parameterized by $ t \in U_{t_0} $, such that:
$$
f_t(h_t(x)) = f_{t_0}(x)
$$
for all $ x $ in a neighborhood of $ 0 \in \bb C^n $ and all $ t \in U_{t_0}$.

The following result follows directly from \cite[Theorem 7.2]{nhan1}.
\begin{thm}\label{thm_criterior_Lipschitztrivial}
Let $F(t,x) = f_t(x): (\bb C^n, 0) \to (\bb C, 0)$ be family of holomorphic function germ where $t$ is in a open connected subset $U \subset \bb C$. If for each $t_0 \in U$ there is continuous vector field $X$ defined on a neighborhood of $(0, t_0)$ in $\bb C^n \times U$ of the form
$$
X(x, t) = \frac{\partial}{\partial t} + \sum_{i=1}^n X_i(x, t) \frac{\partial}{\partial x_i},
$$
and  Lipschitz in $x$ (i.e., there exists a number $C > 0$ with
$$
\|X(x_1, t) - X(x_2, t)\| \leq C \|x_1 - x_2\|
$$
for all $t$), such that $X \cdot F = 0$, then $F$ is a Lipschitz trivial over $U$.
\end{thm}

\subsection{Lipschitz Triviality via Newton's Polyhedron}

Let $x = (x_1, \ldots, x_n) \in \bb C^n$. We denote by $\bar{x} =(\bar{x}_1, \ldots, \bar{x}_n)$ the complex conjugate of $x$. For $\nu = (\nu_1, \ldots, \nu_n) \in  \bb N^n$ we set $x^\nu = x_1^{\nu_1}\ldots x_n^{\nu_n}$  and $\bar{x}^\nu = \bar{x}_1^{\nu_1}\ldots \bar{x}_n^{\nu_n}$.

A function $f: \bb C^n \to \bb C$ is called a \textit{mixed polynomial} if it is of the following form: 

$$f(x, \bar{x}) = \sum_{\nu, \mu} c_{\nu, \mu} x^\nu \bar{x}^\mu, \hspace{0.5 cm} c_{\nu,\mu} \in \bb C^*.$$

Given $w = (w_1,\ldots,w_n) \in \bb N^n$, \textit{the filtration} (with respect to the weight $w$) of a mixed monomial $M = x^\nu \bar{x}^\mu$ is defined by 
$$\text{fil}_w(M) = \sum_{j = 1}^n w_j (\nu_j + \mu_j).$$
The filtration of a mixed polynomial $f$, denoted $\text{fil}_w(f)$ is the minimum of the filtrations of the mixed  monomials appearing in $f$.

A mixed polynomial $f = \sum_{\nu, \mu} c_{\nu, \mu} x^\nu \bar{x}^\mu$ is called radically quasihomogeneous of the type $(w; d)$ if 
$$ \sum_j w_j (\nu_j + \mu_j) = d $$ 
for all $(\nu, \mu)$.

In the case that all $\mu$ are zero, the mixed polynomial $f$ becomes a polynomial in $x$, and the notion of radially quasihomogeneous coincides with the usual notion of quasihomogeneous.

Given a polynomial $f(x)$ and a weight $w$, there is a unique way to express $f$ in the form 
$$ f(x) = H_d(x) + H_{d+1}(x) + \ldots$$
where each $H_k$ (for $k \geq d$) is  a quasihomogeneous polynomial of type $(w;k)$ and $H_d \not \equiv 0$. We call $H_d$ the \textit{initial part of $f$} with respect to $w$. 

The \textit{support} of a mixed polynomial $f = \sum_{\nu, \mu} c_{\nu, \mu} x^\nu \bar{x}^\mu$  is the set  $$\supp (f) = \{(\nu_1 +\mu_1,\ldots, \nu_n + \mu_n)\in \bb Z^n: c_{\nu, \mu} \neq 0\}.$$
We denote by $\Gamma_+(f)$ the \textit{Newton polyhedron} of $f$ which is the convex hull of the set $\bigcup_{\alpha \in \supp(f)} (\alpha + \bb R^n_+)$.  The union $\Gamma(f)$ of the compact faces of $\Gamma_+(f)$ is called \textit{Newton diagram} of $f$. The Newton diagram of $f$ is called \textit{convenient} of the intersection with each coordinate axis is non-empty.  The set of vertices of $\Gamma(f)$ is denoted by $V(\Gamma(f))$.

Let $\sigma$ be a  $(n-1)$-dimensional face of the Newton polyhedron $\Gamma_+(f)$. A \emph{weight} associated to $\sigma$ is a non-zero vector $w = (w_1, \ldots, w_n) \in \mathbb{Q}_{\geq 0}^n$ that is orthogonal to the affine hyperplane containing $\sigma$. the hyperplane has equation equation $\sum_{i = 1}^n w_i x_i = d$, and any mixed monomial $x^\nu \bar{x}^\mu$ where $\nu + \mu$ contained in such a plane is radically quasihomogeneous of the type $(w; d)$. Once the weight $w$ of $\sigma$ is chosen, we often write $\sigma = (w; d)$ and call $d$ the \textit{total weight} of $\sigma$ with respect to $w$. We denote by $w_*$ the maximum value among the $w_i$'s.

We now define the filtration and Newton polyhedron for analytic family of mixed polynomials. Let  
\[
f(t, x, \bar{x}) = \sum_{\nu, \mu} c_{\nu, \mu}(t) \, x^\nu \bar{x}^\mu, \quad \text{with } c_{\nu, \mu}(t) \not\equiv 0,
\]  
where $c_{\nu, \mu}(t)$ is analytic. We regard $ f $ as a family $ f_t $ of mixed polynomials parametrized by $ t $.  

Given a weight $ w \in \mathbb{N}^n $, the filtration of $ f $ with respect to $ w $ is defined as  
\[
\fil_w(f) = \min_{\nu, \mu} \fil_w(x^\nu \bar{x}^\mu).
\]  
Clearly, $ \fil_w(f_t) \geq \fil_w(f) $, and equality holds for generic values of $ t $.  

The support of $ f $ is defined by  
\[
\supp(f) = \left\{ (\nu_1 + \mu_1, \ldots, \nu_n + \mu_n) \in \mathbb{N}^n : c_{\nu, \mu}(t) \not\equiv 0 \right\}.
\]  
The notions of the Newton polyhedron and the Newton diagram of $ f $ are then defined in the same way as in the non-parametric case.

\begin{defn}\rm 
    
\textit{A control function} is a nonzero mixed polynomial $h: \bb C^n \to \bb C$ such that

\begin{enumerate}[label=(\roman*)]
    \item $\Gamma(h)$ is convenient
    \item Near origin, $|h| \sim \rho_{\Gamma (h)}$ where  $\rho_{\Gamma(h)}:= \sum_{\alpha = (\alpha_1, \ldots, \alpha_n)\in V(\Gamma (h))}|x_1|^{\alpha_1}\ldots |x_n|^{\alpha_n}.$
\end{enumerate}

We call  a control function $h$ satisfying $h = \rho_{\Gamma(h)}$ a \textit{standard control function.}
\end{defn}

\begin{lem}\label{lem_control_1}
Let $h$ be a control function. Let $u = (u_1, \ldots, u_n)   \in \Gamma_+(h)$. Then, there are a neighborhood $U$ of $0$ in $\bb C^n$ and  a constant $C >0$ such that for all $x\in U$, 
$$ \rho_{\Gamma(h)}(x)\geq C|x_1|^{u_1}\ldots |x_n|^{u_n}.$$

In particular, if $u \in \Gamma_+(h)\setminus \Gamma(h)$ then 
$$\lim_{x\to 0}\frac{|x_1|^{u_1}\ldots |x_n|^{u_n}}{\rho_{\Gamma(h)} }=0.$$
\end{lem}

\begin{proof} Let $I_u = \{i: u_i \neq 0\}$. Observe that $$|x_1|^{u_1}\ldots |x_n|^{u_n} \neq 0 \Leftrightarrow x_i \neq  0 \text{ for all } i \in I_u.$$

We just need to show that $\rho_{\Gamma(h)} \gtrsim |x_1|^{u_1}\ldots |x_n|^{u_n}$ on the set $X_u: = \{ x_i \neq 0 \text{ for all } i \in I_u\}$.

Now assume in contrast that there exists a real analytic curve $\gamma(t)=\left(\gamma_1(t), \ldots, \gamma_n(t), \gamma_{n+1}(t)\right):[0, \varepsilon) \rightarrow X_u \times \bb R_{\geq 0}$, such that $\gamma(0)=0$, and on $(0, \varepsilon)$, $\gamma_i(t) \neq 0$ for all $i \in I_u$,  $\gamma_{n+1}(t) > 0$, and 
\begin{equation}\label{equ_lem_4.2.1}   
\left|\rho_{\Gamma(h)}(\gamma_1(t), \ldots, \gamma_n(t))\right| < \gamma_{n+1}(t)|\gamma_1(t)|^{u_1}\ldots |\gamma_n(t)|^{u_n},
\end{equation}
Since the right-hand-side of \eqref{equ_lem_4.2.1} is  $> 0$ for all $t \in (0, \varepsilon)$ and $\rho_{\Gamma(h)}$ is continuous, by small pertubation, we may assume $\gamma_i(t) \neq 0$, $i = 1, \ldots, n$ for all $t\in (0, \varepsilon)$. 

Thus, we can assume that  $\left|\gamma_1(t)\right| \sim t^{a_1}, \ldots,\left|\gamma_n(t)\right| \sim t^{a_n},\left|\gamma_{n+1}(t)\right| \sim t^b$ and
where $a_i >0$ for $i = 1, \ldots, n$, $b>0$. Then, \eqref{equ_lem_4.2.1} induces that  
$$
\sum_{\alpha \in \Gamma(h)} t^{a_1 \alpha_1} \ldots t^{a_n \alpha_n} \lesssim t^b t^{a_1 u_1} \ldots t^{a_n u_n} .
$$

As $t$ small enough, we have
$$
\sum_{\alpha \in \Gamma(h)} t^{\langle a, \alpha\rangle}<t^{\langle a, u\rangle} .
$$

This implies 
\begin{equation}\label{equ_4.1.0}
    \langle a, u\rangle<\inf _{\alpha \in V(\Gamma(h))}\langle a, \alpha\rangle
\end{equation}

Fix $a$ and consider the linear function  $\langle a, v\rangle$ where $v\in \Gamma_+(h)$. This function reaches a minimum at one of the vertices of $\Gamma(h)$, and hence $$ \langle a, u\rangle \geq \inf_{v \in \Gamma_+(h)}  \langle a, v\rangle =  \inf _{v\in V(\Gamma(h))} \langle a, v\rangle.$$ This contradicts \eqref{equ_4.1.0}. Therefore, the first part of the lemma is proved. 

Now let us prove the second part.  Assume $u \in \Gamma_+(h) \setminus \Gamma(h)$.  Then, there is a nonzero vector $v = (v_1, \ldots, v_i)$ with $v_i\geq 0$ sufficiently small such that  $(u - v)\in \Gamma_+(h)$. This implies that
\begin{align*}
    \frac{|x_1|^{u_1}\ldots |x_n|^{u_n}}{\rho_{\Gamma(h)}(x)} & = \frac{(|x_1|^{u_1-v_1}\ldots |x_n|^{u_n-v_n}) (|x_1|^{v_1}\ldots |x_n|^{v_n})}{\rho_{\Gamma(h)}(x)}\\
    & \lesssim|x_1|^{v_1}\ldots |x_n|^{v_n} \to 0 \text{ as } x\to 0,
\end{align*}
where $C>0$ is a some constant. 
\end{proof}

We have the following result:

\begin{thm}\label{thm_checking_triviality_1}
Let $\Delta \subset \bb C$ be an open subset and let  $F_t: (\bb C^n, 0)\to (\bb C, 0)$, $F_t(x) = f(x) + t g(x)$, $t \in \Delta$ be a family of germs of polynomials with isolated singularity, such that all $F_t$ share the same Newton diagram. Let $h_t(x) = h(t,x)$ be an analytic family of control functions of the form
\begin{equation}\label{equ_4.1.2}
    h(t,x) = \sum_{i = 1}^n P_i(t,x) \frac{\partial F_t}{\partial x_i}(t,x),
\end{equation}
where each $P_i$ is an analytic family of germs of mixed polynomials, and the Newton diagram of $h_t$ is independent of $t$.

For every $(n-1)$-dimensional compact face $\sigma$ of $\Gamma(h)$ and for every $i = 1, \ldots, n$,
    \begin{equation}\label{equ_4.1.3}
        \fil_{w_\sigma} (g) + \fil_{w_\sigma} (P_{i})  -  \fil_{w_\sigma}(h)  - w_{\sigma,*} \geq 0,
    \end{equation}
    where $w_\sigma$ is a weight associate to $\sigma$,
    then $F_t$ is Lipschitz trivial over $\Delta$.
\end{thm}

\begin{proof}
It is obvious that 
\begin{align}\label{eq_lem_2}
h\frac{\partial F_t}{\partial t} =  \sum_i^n \left(P_i(t,x) \frac{\partial F_t}{\partial t}  \right)\frac{\partial F_t}{\partial x_i} .
\end{align}

Put 
$$A_i(t,x) = 
\begin{cases}
\frac{P_i(t,x)  \frac{\partial F_t}{\partial t}}{h(t,x)},& \text{ if } x \neq 0\\
0, & \text{if } x = 0 
\end{cases}
$$

Set
$$X(t, x) = -\frac{\partial}{\partial t}+\sum_{i=1}^n A_i(t,x)\frac{\partial}{\partial x_i}.$$

Then, 
$$ X.f = 0.$$

By Theorem \cite[Theorem 7.1]{nhan1}, to prove that $ F_t $ is Lipschitz trivial, it suffices to show that, on neighborhood of 
$\Delta\times \{0\}$, $X$ is continuous in $(t,x)$ and locally Lipschitz in $x$.

\textbf{Claim 1:} $X$ is continuous on a neighborhood of $\Delta \times \{0\}$

We write $P_i(t,x)g(x)$ as the form  
$$P_i(t,x)g(x) = \sum_{(\nu, \mu)\in I_i}c_{\nu, \mu} (t) x^\nu \bar{x}^\mu. $$ 

Fix $t_0\in \Delta$,  take small neighborhood $\Delta_{t_0}$ of $t_0$ in $\bb C$ such that  $\overline{\Delta}_{t_0}  \subset U$. 
Set $$\delta_{t_0} := \max_{(\nu, \mu)\in I_i,t\in\overline{\Delta}_{t_0}}  c_{\nu, \mu}(t).$$ Note that such a $\delta_{t_0}$ exists since $c_{\nu, \mu}(t)$ are continuous functions. 

It is clear that for every $x\in \bb C^n$ and every $t\in \Delta_{t_0}$, 
$$|P_i(t,x) g(x)| \leq \delta_{t_0} \sum_{(\nu, \mu)\in I_i} |x|^\nu |\bar{x}|^\mu =  \delta_{t_0} \sum_{(\nu, \mu)\in I_i} |x|^{\nu+\mu}.$$

Since $h_t$ are control functions and $\Gamma(h_t)$ is independent of $t$, for each $t$ there is a neighborhood of the origin in $\bb C^n$ such that 
$$\rho_{\Gamma(h)} = \rho_{\Gamma(h_t)} \sim |h_t(x)|.$$ Note that the constant for the relation $\sim$ above is depending on $t$. However, because $h(t,x)$ is continuous, if we fix $t = t_0$, shrinking $\Delta_{t_0}$ smaller if necessary, it is possible to choose an open neighborhood $U_{t_0}$ of $0$ in $\bb C^n$ such that on $ \Delta_{t_0} \times U_{t_0}$
$$ \rho_{\Gamma(h)} (x)\sim |h(t,x)|.$$

Note that $g(x) = \frac{\partial F_t}{\partial t} (t,x)$ and by the assumption  \eqref{equ_4.1.3} $$\fil_{w_\sigma}(P_{i} g) = \fil_{w_\sigma} (g) + \fil_{w_\sigma} (P_{i}) > \fil_{w_\sigma} h$$ for every $(n-1)$-dimensional compact face $\sigma$ of $\Gamma_+(h)$. Together with the convenience of $\Gamma(h)$, this implies that $\supp (P_{i} g) \subset \Gamma_+(h)\setminus \Gamma(h)$. Hence,   $\mu+\nu \subset \Gamma_+(h)\setminus \Gamma(h)$ for all $(\nu, \mu) \in I_i$. By Lemma \ref{lem_control_1}, 
$$\lim_{x\to 0}\frac{|x|^{\nu +  \mu}}{\rho_{\Gamma(h)}(x)} =  0.$$

 On $\Delta_{t_0}\times U_{t_0} $ we have $\rho(x) \sim |h(t,x)|$, and $h^{-1}(0) =  \Delta_{t_0} \times \{0\}$. Therefore,

 \begin{align*}
\lim_{(t,x) \to \{0\} \times \Delta_{t_0}}|A_{i}(t,x)| & = \lim_{(t,x) \to \{0\} \times \Delta_{t_0}} \frac{|P_i(x)g (x)|}{|h(t,x)|} \\
& \lesssim \lim_{(t,x) \to \{0\} \times \Delta_{t_0}} \frac{\sum_{\nu, \mu} |x|^{\nu +  \mu}}{\rho_{\Gamma(h)}(x)} = 0
 \end{align*} 

Thus, $X$ is continuous.   

\textbf{Claim 2:} $X$ is Lipschitz in $x$ on a neighborhood of $\Delta \times \{0\}$.

We will show that with a $t_0\in \Delta$ fixed, on the neighborhood $\Delta_{t_0}\times U_{t_0}$ as in the proof of Claim 1, $X$ is Lipschitz in $x$. Since $h_t^{-1}(0) = \{0\}$, to show that $X$ is Lipschitz in $x$ on $\Delta_{t_0}\times U_{t_0}$, we just need to show that all first derivatives of $X$ in $x_j$ and $\bar{x}_j$ are bounded by a constant independent of $t \in \Delta_{t_0}$.

For $(t,x)\not \in h^{-1}(0)$, computation gives
$$
\frac{\partial X}{\partial x_j} =  \frac{\partial A_{i}}{\partial x_j}
= \frac{
\left( \dfrac{\partial P_{i}}{\partial x_j} \cdot g + P_{i} \cdot \dfrac{\partial g}{\partial x_j} \right) \cdot h
- P_{i} \cdot g \cdot \dfrac{\partial h}{\partial x_j}
}{
h^2
}
$$

Note that for every  $(n-1)$-dimensional compact face $\sigma$ of $\Gamma(h)$,
$\fil_{w_\sigma}(h^2) = 2 \fil_{w_\sigma}(h)$. In addition, since $h_t$ is a control function, so is $h_t^2$.

We have  
\begin{align*}
 \fil_{w_\sigma}  \left( \left( \dfrac{\partial P_{i}}{\partial x_j} \cdot g + P_{i} \cdot \dfrac{\partial g}{\partial x_j} \right) \cdot h
- P_{i} \cdot g \cdot \dfrac{\partial h}{\partial x_j}\right) & \geq \fil_{w_\sigma}(P_{i}) + \fil_{w_\sigma}(g) + \fil_{w_\sigma}(h) - w_{\sigma,j} \\
& \overset{\eqref{equ_4.1.3}}{\geq} \fil_{w_\sigma}(h^2).
\end{align*}

Similar arguments as in the proof of Claim 1 and by Lemma \ref{lem_control_1}, taking $U_{t_0}$ and $\Delta_{t_0}$ smaller if necessary, we conclude that $\frac{\partial X}{\partial x_j}$ is bounded on $\Delta_{t_0}\times U_{t_0} \setminus h^{-1}(0)$. We can show the boundedness of  $\frac{\partial X}{\partial \bar{x}_j}$ similarly. This completes the proof. 
\end{proof}

The following is a sufficient condition for a function $h$ to be a control function. 

\begin{lem}[{\cite[Lemma 7.5]{nhan1}}]\label{lem_control0} Let $h: (\bb C^n, 0) \to (\bb C, 0)$ be the germ of a mixed polynomial such that  $\Gamma(h)$ is convenient and for every compact face $\sigma$ of $\Gamma(h)$, the equation $h|_\sigma(x) = 0$, near the origin, has no solution in $(\bb C^*)^n$.  Then, there is a constant $C>0$ such that in a neighborhood of the origin we have 
$$ |h(x)| \geq C \rho_{\Gamma(h)}(x).$$
\end{lem}

It directly follows from Lemma \ref{lem_control_1} and  Lemma \ref{lem_control0} that 

\begin{lem}\label{lem_control2}
    Let $h$ be given as in Lemma \ref{lem_control0}. Then $h$ is a control function. 
\end{lem}

We also require the following result for the proof of Theorem~\ref{thm_checking_triviality_main}.

\begin{thm}[{\cite[Proposition 7.1]{nsr}}]\label{thm_checking_trivialiy_2}
Let $g: (\bb C^n,0) \to (\bb C, 0)$ be a germ of a quasihomogeneous polynomial of type $(w_1, \ldots, w_n; d)$. Let $U$ be an open connected subset of $\bb C$ and let $F(t,x) = F_t(x)  = g(x) + t\theta(x)$, $t\in U$ be an analytic family of germs of polynomials. Suppose that the initial part of $F_t$ with respect to $(w_1, \ldots, w_n)$ is has an isolated singularity for every $t$. If $\mathrm{fil}(\theta) \geqslant d + \max_{i,j} \{ w_i - w_j \}$, then $F$ is Lipschitz trivial over $U$.
\end{thm}

\begin{thm}\label{thm_checking_triviality_main}
    The following families are Lipschitz trivial with respect to the direction of $t$.
     \begin{enumerate}
    \item $Z_{1,0}: x^3 y + sx^2y^3 + txy^6 + y^7$, $4s^3 + 27 \neq 0$
    \item $Z_{1, 1}: x^3 y + x^2 y^3 + sy^8 + t y^9$, $s \neq 0$ 
     \item $Z_{17}: x^3y + y^8 + sxy^6 + txy^7$, $s\neq 0$
    
   \item  $W_{1,0}: x^4 +  sx^2y^3 + t x^2 y^4  +  y^{6}$, $s^2\neq 4$ 
   \item $W_{1,1}: x^4 + x^2y^3  +  sy^{7} + t y^8$, $s \neq 0$
   \item $W_{1,2}: x^4 + x^2y^3  +  sy^{8} + t y^9$, $s \neq 0$
   \item $W_{17}: x^4 + xy^5 + sy^7 + ty^8$, $s\neq 0$
    \item $W_{18}: x^4 + y^7 + sx^2y^4 + tx^2 y^5$, $s\neq 0$
    
    \item $W_{1, 2q-1}^{\#}: (x^2 + y^3)^2 + s x y^{4+q} + t xy^{5 + q}$, $s \neq 0$, $q>0$
    \item $W_{1, 2q}^{\#}: (x^2 + y^3)^2 + sx^2 y^{3+q} + t x^2 y^{4 + q}$, $s \neq 0$, $q >0$
    \end{enumerate}
\end{thm}

\begin{proof} 
\textit{Proof of (1)}:  Consider the family
$$f_{s,t}(x,y) = x^3 y + sx^2y^3 + txy^6 + y^7, \ 4s^3 + 27 \neq 0.$$
Put $g_{s}(x, y) = x^3 y + sx^2y^3 + y^7$ and $\theta(x,y) = xy^6$. Let $w = (2, 1)$ be a weight. 
It is clear that $g_s$ is a quasihomogeneous germ of the type $(2, 1; 7)$ and $\fil_w (\theta) = 8$. Applying Theorem \ref{thm_checking_trivialiy_2}, it follows that $f_{s,t} = g_{s} + t h$ is Lipschitz trivial with respect to the parameter $t$. 

Theorem~\ref{thm_checking_trivialiy_2} does not apply to families (2)--(10). 
Instead, for these cases we use Theorem~\ref{thm_checking_triviality_1}. 
The idea is as follows: for a fixed value of $s$, we construct a family of control functions $h(s,t,x,y)$ satisfying the conditions of Theorem~\ref{thm_checking_triviality_1}. 
In most cases, the family of control functions takes the form
\begin{equation}\label{equ_lip_trivial_1}
    h(s,t,x,y) = u_1(s,t,x,y)\,|x|^{2a} + u_2(s,t,x,y)\,|y|^{2b},
\end{equation}
with $a \leq b$, where
\begin{enumerate}[label=(\roman*)]
  \item $u_1 x^a = \alpha_1 \frac{\partial f}{\partial x} + \beta_1 \frac{\partial f}{\partial y}$ and 
        $u_2 y^b = \alpha_2 \frac{\partial f}{\partial x} + \beta_2 \frac{\partial f}{\partial y}$,  
        where for $i=1,2$, $u_i(s,t,x,y)$ are units for all $s,t$, and $\alpha_i,\beta_i$ are polynomials;
  \item $\min_i\{\fil_w \alpha_i, \fil_w \beta_i\} + \fil_w (\partial f/\partial t) - b/a - b \geq 0$, 
        where $w=(b/a,1)$. 
\end{enumerate}
Condition (i) can be verified with the help of \textsc{Singular}.  

We now explain why such a family of control functions satisfies the assumptions of Theorem~\ref{thm_checking_triviality_1}. 
Indeed,
\begin{align*}
    h(s,t,x,y) 
    &= \bar{x}^a u_1 x^a + \bar{y}^b u_2 y^b \\
    &= (\bar{x}^a \alpha_1 + \bar{y}^b \beta_1)\,\frac{\partial f}{\partial x} 
       + (\bar{x}^a \alpha_2 + \bar{y}^b \beta_2)\,\frac{\partial f}{\partial y} \\
    &= P_1(s,t,x,y)\,\frac{\partial f}{\partial x} + P_2(s,t,x,y)\,\frac{\partial f}{\partial y}.
\end{align*}
The Newton diagram $\Gamma(h)$ has a unique compact one-dimensional face with weight $w=(b/a,1)$. 
Note that $\fil_w(h)=2b$ and $w_*=b/a$. 
By condition (ii), we have 
\[
\fil_w P_j\geq  b + \min_i \{\fil_w \alpha_i, \fil_w \beta_i\}, j = 1, 2.
\]
It follows that 
$$\fil_w (P_i) + \fil_w(\partial f /\partial t) - \fil_w h - w_* \geq b + \min_i \{\fil_w \alpha_i, \fil_w \beta_i\} + \fil_w(\partial f /\partial t) - 2b - b/a \overset{\text{(ii)}}{\geq}  0$$
and hence the requirement \eqref{equ_4.1.3} is satisfied. 
Therefore, $f$ is bi-Lipschitz trivial in $t$.  

We now give a detailed proof of (4).
$$W_{1,0}: f(s,t,x,y) = x^4 + s x^2 y^3 + t x^2 y^4 + y^6, \quad s^2 \neq 4.$$

Fix $s$ with $s^2 \neq 4$.

\textbf{Case 1: $s=0$.} In this case, 
$$f(s,t,x,y) = x^4 + t x^2 y^4 + y^6.$$
This family is quasihomogeneous of type $(3,2;12)$. 
By Theorem~\ref{thm_checking_trivialiy_2}, we conclude that $f$ is bi-Lipschitz trivial.  

\textbf{Case 2: $s \neq 0$.} We compute
\[
\frac{\partial f}{\partial x} = 4x^3 + 2s x y^3 + 2t x y^4,
\qquad 
\frac{\partial f}{\partial y} = 3s x^2 y^2 + 4t x^2 y^3 + 6y^5.
\]
We will construct a family of control functions $h(s,t,x,y)$ as in \eqref{equ_lip_trivial_1} with $a=5$ and $b=8$.  
To express $x^5$ and $y^8$ in the form (i), one may use \textsc{Singular}. 

SINGULAR computation: 
\begin{verbatim}
ring R = (0, s,t), (x, y), ds;
poly f = x4 + sx2y3 + tx2y4 + y6;
ideal I = jacob(f); 
division(x5, I);
[1]:
   _[1,1]=-1/(2519424s7-10077696s5)*x3
          +(-7t)/(7558272s8-60466176s6+120932352s4)*x3y
          -1/(1259712s8-10077696s6+20155392s4)*xy3 
          +(-t2)/(1889568s9-15116544s7+30233088s5)*x3y2
          +(-t)/(1259712s9 -10077696s7+20155392s5)*xy4
   _[2,1]=1/(3779136s7-30233088s5+60466176s3)*x2y
          +(t)/(1889568s8-15116544s6+30233088s4)*x2y2
          +(t2)/(3779136s9-30233088s7+60466176s5)*x2y3
[2]:
   _[1]=0
[3]:
   _[1,1]=-1/(629856s7-2519424s5)+(-7t)/(1889568s8-15116544s6+30233088s4)*y
          +(-t2)/(472392s9-3779136s7+7558272s5)*y2
> division(y8, I);
[1]:
   _[1,1]=(s)/(108s4-864s2+1728)*xy2+(t)/(81s4-648s2+1296)*xy3
   _[2,1]=-1/(81s4-648s2+1296)*x2-1/(162s3-648s)*y3
          +(-3s2t-16t)/(486s6-3888s4+7776s2)*y4+(32t2)/(729s7-5832s5+11664s3)*y5
[2]:
   _[1]=0
[3]:
   _[1,1]=-1/(27s3-108s)+(-3s2t-16t)/(81s6-648s4+1296s2)*y
          +(128t3)/(729s7-5832s5+11664s3)*x2+(64t2)/(243s7-1944s5+3888s3)*y2
> 
\end{verbatim}
It follows from the computations above that $x^5$ and $y^8$ can be written in the form (i). 
The monomials appearing in $\alpha_1, \beta_1, \alpha_2, \beta_2$ are, respectively:
\[
\{x^3,\, x^3y,\, xy^3,\, x^3y^2,\, xy^4\}, \quad 
\{x^2y,\, x^2y^2,\, x^2y^3\}, \quad
\{xy^2,\, xy^3\}, \quad
\{x^2,\, y^3,\, y^4,\, y^5\}.
\]

The monomials of $\alpha_i$ and $\beta_i$ can also be extracted directly in \textsc{Singular}. 
For example, the following code lists the monomials in $\alpha_1$:
\begin{verbatim}
list L = division(x5, I); 
poly P = matrix(L[1])[1,1];
matrix M = coef(P, xy); 
for (int i = 1; i <= ncols(M); i++) {
    print(M[1,i]); 
};
xy4
x3y2
xy3
x3y
x3
\end{verbatim}

With the weight $w=(8/5,1)$ we compute
\[
\min_{i}\{\fil_w(\alpha_i), \fil_w(\beta_i)\} = \fil_w(y^3) = 3,
\qquad 
\fil_w\!\left(\frac{\partial f}{\partial t}\right) 
   = \fil_w(x^2y^4) = \tfrac{36}{5}.
\]
Hence,
\[
3 + \frac{36}{5} - \frac{8}{5} - 8 \;=\; \frac{3}{5} > 0,
\]
so condition (ii) is satisfied. 
Therefore, $f$ is bi-Lipschitz trivial in $t$.

\textit{Proof of (2), (3), (5), (6), and (8).}  
For these families, we take control functions $h(s,t,x,y)$ of the form \eqref{equ_lip_trivial_1} with $a=5$ and $b=10$. 
In this setting, $w=(2,1)$ is the weight associated with the unique one-dimensional compact face of $\Gamma(h)$. 
Note that the functions $u_i, \alpha_i, \beta_i$ $(i=1,2)$ depend on each specific case.  
The monomials appearing in $\alpha_i$ and $\beta_i$, together with $\min_{i}\{\fil_w(\alpha_i), \fil_w(\beta_i)\}$ and $\fil_w(\partial f/\partial t)$ for each case, are listed in Table~\ref{tab1:list of monomials}. In the list of monomials, the symbol ``$\ldots$'' indicates that additional monomials are generated from one of the preceding ones.

\begin{center}
\footnotesize
\begin{table}[h!]
\caption{List of monomials in $\alpha_i$ and $\beta_i$, $i = 1, 2$; $\min_{i}\{\fil_w(\alpha_i), \fil_w(\beta_i)\}$ and $\fil_w(\partial f/\partial t)$ for families (2), (3), (5), (6) and (8)}
\begin{tabular}{|c|c|c|c|c|}
\hline
\textbf{Family} & \textbf{$\alpha_i, \beta_i$} & \textbf{Monomials in $\alpha_i, \beta_i$} & $\min_{i}\{\fil_w(\alpha_i), \fil_w(\beta_i)\}$ &  \textbf{$\fil_w \frac{\partial f}{\partial t}$}\\
\hline
\multirow{3}{*}{(2): $Z_{1,1}$} & $\alpha_1$ & $x^2y, xy^4, y^6, \ldots$ & \multirow{3}{*}{3}  &\multirow{3}{*}{9} \\
\cline{2-3} 
                  & $\beta_1$ &  $x^2, xy^2, \ldots$ &  & \\
\cline{2-3}
                  & $\alpha_2$ &  $ x^2,\ xy^2,y^5, \ldots $  & &\\
\cline{2-3}
                  & $\beta_2$ & $xy, y^3,\ldots$  &  & \\
\hline
\multirow{3}{*}{(3): $Z_{17}$} & $\alpha_1$ & $ x^3, xy^4, x^2y^2, y^6, \ldots $ & \multirow{3}{*}{3} & \multirow{3}{*}{9} \\
\cline{2-3} 
                  & $\beta_1$ &  $x^2, xy^3, y^5, \ldots$  &   & \\
\cline{2-3}
                  & $\alpha_2$ &  $ x^2, xy^2, y^5, \ldots $   &    &\\
\cline{2-3}
                  & $\beta_2$ & $xy, y^3, \ldots$  &    &\\
\hline
\multirow{3}{*}{(5): $W_{1,1}$} & $\alpha_1$ & $x^2, y^4, \ldots$ & \multirow{3}{*}{3} &\multirow{3}{*}{9} \\
\cline{2-3} 
                  & $\beta_1$ &  \shortstack{$xy$}  &   & \\
\cline{2-3}
                  & $\alpha_2$ &  $ xy^3, \ldots$   &  & \\
\cline{2-3}
                  & $\beta_2$ & $x^2y, y^4,\ldots$  &  & \\
\hline
\multirow{3}{*}{(6): $W_{1,2}$} & $\alpha_1$ & $x^2, y^5, \ldots$ & \multirow{3}{*}{3}  &\multirow{3}{*}{9} \\
\cline{2-3} 
                  & $\beta_1$ &  $xy$  &   & \\
\cline{2-3}
                  & $\alpha_2$ &  $ xy^2, \ldots$   &  & \\
\cline{2-3}
                  & $\beta_2$ & $x^2, y^3, \ldots$  &  & \\
\hline
\multirow{3}{*}{(8): $W_{18}$} & $\alpha_1$ & $x^2, y^4, \ldots$ & \multirow{3}{*}{4}  &\multirow{3}{*}{9} \\
\cline{2-3} 
                  & $\beta_1$ &  $xy^2, \ldots$  &  & \\
\cline{2-3}
                  & $\alpha_2$ &  $xy^4, \ldots$    & & \\
\cline{2-3}
                  & $\beta_2$ & $x^2y, y^4, \ldots$  &  & \\
\hline
\end{tabular}

\label{tab1:list of monomials}
\end{table}
\end{center}
It is straightforward to verify that condition (ii) of \eqref{equ_lip_trivial_1} is satisfied.

Several other cases can be proved by the same method. 
For example:  
(7) $W_{17}$, where we take control functions of the form \eqref{equ_lip_trivial_1} with $a=10$, $b=10$;  
(9) $W_{1,1}^{\#}$ with $a=6$, $b=8$;  
(10) $W_{1,2}^{\#}$ with $a=6$, $b=9$, and $W_{1,4}^{\#}$ with $a=7$, $b=10$.

There remain only two cases: (9) $W_{1,2q-1}^{\#}$ with $q\geq 2$, and (10) $W_{1,2q}^{\#}$ with $q\geq 3$. 
Since these cases consist of finitely many families of germs, it is not feasible to use \textsc{Singular} to check them all. 
Therefore, we handle these cases by hand.

\textit{Proof of (9) with $q\geq 2$}.

Consider the family: 
$$W_{1,2q-1}^{\#}: f(s,t, x, y) =  (x^2 + y^3)^2 + s x y^{q+4} + t xy^{q+5}, \ s \neq 0, q\geq 2$$

\begin{lem}\label{lem_(6)} There are units $u_i$ and polynomials $a_i$, $b_i$, $i = 1, 2$ such that 

(i) $u_1 u^{q+7} = \alpha_1 \partial f/\partial x + \beta_1 \partial f/\partial y$

(ii)  $u_2 x^{q+4} = \alpha_2 \partial f/\partial x + \beta_2 \partial f/\partial y$.

(iii) In particular, the monomials appearing in all $\alpha_i$ and $\beta_i$, $i = 1, 2, 3$ includes $x^2, xy, y^3$ and other monomials generated from these.
\end{lem}
Consider the family of control functions: 
$$ h(s,t, x, y) = u_1 |x|^{2q+8} + u_2 |y|^{2q+14}.$$
It is also of the form \eqref{equ_lip_trivial_1} with $a = q+4$ and $b = q+7$. 
And, $w = (1 + \frac{3}{q+4}, 1)$ the weight of the unique compact $1$-dimensional face in $\Gamma(h)$. It is clear that 
$$\min_{i}\{\fil_w \alpha_i, \fil_w \beta_i\} =  \fil_w xy =2 + \frac{3}{q+4}.$$
In addition, $\fil_w \partial f/\partial t = \fil_w xy^{q+5} = q+6 + \frac{3}{q+4}$. We have
\begin{align*}
    \min_{i}\{\fil_w \alpha_i, \fil_w \beta_i\} + \fil_w \partial f/\partial t  - b/a - b  & =  2 + \frac{3}{q+4} +  q+6 + \frac{3}{q+4} - (1 + \frac{3}{q+4}) - (7+q) \\
    & = \frac{3}{q+4} \geq 0
\end{align*}
Then, (ii) in \eqref{equ_lip_trivial_1} holds, and hence $f$ is bi-Lipschitz trivial in $t$.
This ends the proof of $(9)$. 

Let us now back with the proof of Lemma \ref{lem_(6)}. 

\begin{proof}[Proof of Lemma \ref{lem_(6)}]
(i) We first compute the partial derivatives:
\[
\begin{aligned}
\frac{\partial f}{\partial x} &= 4x^3 + 4x y^3 + s\,y^{4+q} + t\,y^{5+q},\\[4pt]
\frac{\partial f}{\partial y} &= 6x^2 y^2 + 6y^5 + s(4+q)\,x\,y^{3+q} + t(5+q)\,x\,y^{4+q}.
\end{aligned}
\]

Consider the linear combination
\[
A := y^3 \frac{\partial f}{\partial x} - \tfrac{2}{3}x y \frac{\partial f}{\partial y}.
\]
Expanding and canceling terms ($4x^3y^3$ and $4xy^6$), we obtain
\begin{equation}\label{eq_exp_2}
A = s y^{q+7} + t y^{q+8}
- \tfrac{2s(4+q)}{3} x^{2} y^{4+q}
- \tfrac{2t(5+q)}{3} x^{2} y^{5+q}.
\end{equation}
From the expression for $\frac{\partial f}{\partial y}$ we have
\begin{equation}\label{eq_exp_3}
6x^2y^2 = \frac{\partial f}{\partial y} - 6y^5 - s(4+q)\,x\,y^{3+q} - t(5+q)\,x\,y^{4+q}.
\end{equation}
Multiplying by $y^{2+q}/6$ yields
\begin{equation}\label{eq_exp_4}
    x^2 y^{4+q}
= \tfrac{1}{6}y^{2+q}\frac{\partial f}{\partial y} - y^{q+7}
- \tfrac{s(4+q)}{6}x\,y^{2q+5}
- \tfrac{t(5+q)}{6}x\,y^{2q+6}.
\end{equation}

Form \eqref{eq_exp_2} and \eqref{eq_exp_4}, we have 
\begin{equation}
    A = -\frac{s(4+q)}{9} y^{2+q} \frac{\partial f}{\partial y} + (s(1+\frac{2(4+q)}{3}) + \hot) y^{q+7}
\end{equation}
It follows that 
$$ y^3 \frac{\partial f}{\partial x} + \left(s\left(1+\frac{2(4+q)}{3}y^{q+2}\right) -\frac{2}{3} xy \right) \frac{\partial f}{\partial y} =  \left(s\left(1+\frac{2(4+q)}{3}\right) + \hot \right) y^{q+7}.$$
Thus, $y^{q+7} \in \jac(f)$.

\medskip
\noindent\emph{Claim 1.} $x^3y^{4+q}\in \jac(f)$ for $q\geq 2$, and $x^{q-2}y^9\in \jac(f)$ for $q\geq 4$.

\smallskip
Consider a monomial $x^n y^m$ with $m,n \geq 2$. We have
\[
x^n y^m = \tfrac{1}{6} x^{n-2} y^{m-2} \frac{\partial f}{\partial y} - x^{n-2} y^{m+3} 
- \tfrac{(4+q) s}{6} x^{n-1} y^{m+1+q} 
- \tfrac{(5+q) t}{6} x^{n-1} y^{m+2+q}.
\]
Thus, for $q\geq 2$,
\[
x^ny^m \in \langle \frac{\partial f}{\partial y}, x^{n-2} y^{m+3}\rangle.
\]
Replacing $(n,m)$ by $(n-2,m+3)$ repeatedly (as long as $n-2k\geq 0$), we obtain
\[
x^ny^m \in \langle \frac{\partial f}{\partial y}, x^{n-2k} y^{m+3k}\rangle,\quad k=1,2,\dots.
\]
Since (i) shows $y^{q+7}\in \jac(f)$, it follows that if $m+3k\geq q+7$ then $x^{n-2k}y^{m+3k}\in \jac(f)$ and hence $x^ny^m\in \jac(f)$.
\begin{itemize}
    \item If $n=2k$, the condition is
\[
m+3k \geq q+7 \quad\Leftrightarrow\quad \tfrac{2m+3n-14}{2}\geq q.
\]
\item  If $n=2k+1$, the condition is
\[
m+3k \geq q+7 \quad\Leftrightarrow\quad \tfrac{2m+3n-17}{2}\geq q.
\]
\end{itemize}

For $x^3y^{4+q}$ ($n=3$, $m=4+q$), this inequality holds with equality, so $x^3y^{4+q}\in \jac(f)$. Similarly, $x^{q-2}y^9\in \jac(f)$ for $q\geq 4$. This proves Claim~1.

\medskip
\noindent\emph{Claim 2.} $xy^{6+q}\in \jac(f)$.

\smallskip
Indeed, consider
\begin{align*}
\tfrac{1}{s} x y^2 \frac{\partial f}{\partial x} - \tfrac{2}{3s} x^2 \frac{\partial f}{\partial y} + \tfrac{4+q}{6} y^{3+q} \frac{\partial f}{\partial x} - \tfrac{2}{3s} x^2 y \frac{\partial f}{\partial y}
&= -\tfrac{2((4+q)s + (5+q) t)}{3s} x^3 y^{4+q} - \tfrac{2(5+q) t}{3s} x^3 y^{5+q} \\
&\quad+ \tfrac{11 + 2q}{3} x y^{6+q} + \tfrac{t}{s} x y^{7+q} \\
&\quad+ \tfrac{(4+q) s}{6} y^{7+2q} + \tfrac{(4+q) t}{6} y^{8+2q}.
\end{align*}
Since $x^3y^{4+q}, y^{7+q}\in \jac(f)$, it follows that $xy^{6+q}\in \jac(f)$. This proves Claim~2.

\medskip
\noindent(ii) We now prove $x^{q+4}\in \jac(f)$.

Recall the partial derivatives
\[
\begin{aligned}
\frac{\partial f}{\partial x} &= 4x^3 + 4x y^3 + s\,y^{4+q} + t\,y^{5+q},\\[4pt]
\frac{\partial f}{\partial y} &= 6x^2 y^2 + 6y^5 + s(4+q)\,x\,y^{3+q} + t(5+q)\,x\,y^{4+q}.
\end{aligned}
\]

Consider
\[
B := x^{1+q}\frac{\partial f}{\partial x} - \tfrac{2}{3}x^{q}y \frac{\partial f}{\partial y}.
\]
Expanding, the $x^{q+2}y^3$ terms cancel, yielding
\[
B = 4x^{q+4} - 4x^{q}y^{6} - \tfrac{2q+5}{3}s\,x^{1+q}y^{4+q} - \tfrac{2q+7}{3}t\,x^{1+q}y^{5+q}.
\]
Thus
\[
4x^{q+4} = B + R,
\]
where
\[
R := 4x^{q}y^6
+ \tfrac{2(4+q)s}{3}x^{1+q}y^{4+q} + \tfrac{2(5+q)t}{3}x^{1+q}y^{5+q}.
\]

Since $B\in \jac(f)$, it suffices to show $R\in \jac(f)$. For $q\geq 2$, Claim~1 implies $x^{1+q}y^{4+q}\in \jac(f)$. Hence it remains to show $x^qy^6\in \jac(f)$.

From \eqref{eq_exp_3} we obtain
\[
x^{q} y^{6} = \tfrac{x^{q-2} y^{4}}{6} f_{y} - x^{q-2} y^{9} - \tfrac{s(4+q)}{6} x^{q-1} y^{7+q} - \tfrac{t(5+q)}{6} x^{q-1} y^{8+q}.
\]
The right-hand side belongs to $\langle \frac{\partial f}{\partial y}, x^{q-2}y^9, y^{q+7}\rangle$. Since $y^{q+7}\in \jac(f)$, it follows that $x^qy^6 \in \langle \frac{\partial f}{\partial x},\frac{\partial f}{\partial y}, x^{q-2}y^9\rangle$.

Finally, $x^{q-2}y^9\in \jac(f)$ for all $q\geq 2$: indeed, for $q=2$ we have $x^{0}y^9=y^9=y^{q+7}\in \jac(f)$; for $q=3$, $x^{1}y^9=xy^{q+6}\in \jac(f)$ by Claim~2; and for $q\geq 4$, $x^{q-2}y^9\in \jac(f)$ by Claim~1. 

Thus $R\in \jac(f)$, and therefore $x^{q+4}\in \jac(f)$. 

Finally, (iii) follows from the constructions of the proofs of (i) and (ii).  The lemma is proved.

\textit{Proof of (10) with $q >2$:}

Consider the family
$$W_{1,2q}^{\#}: f(s,t, x, y) =  (x^2 + y^3)^2 + s x^2 y^{q+3} + t x^2y^{q+4}, \ s \neq 0, q\geq 3$$

\begin{lem}\label{lem_(7)}
There are units $u_i$ and polynomials $a_i$, $b_i$, $i = 1, 2, 3$ such that 

(i) $u_1 y^{q+8} = \alpha_1 \partial f/\partial x + \beta_1 \partial f/\partial y$

(ii)  $u_2 xy^{q+6} = \alpha_2 \partial f/\partial x + \beta_2 \partial f/\partial y$.

(iii)  $u_3 x^{q+4} = \alpha_3 \partial f/\partial x +\beta_3 \partial f/\partial y$.

(iv) Monomials appearing in all $\alpha _i$ and $\beta_i$, $i = 1, 2, 3$ contain $x^2, xy, y^3$ and other monomials generated from these.
\end{lem}

We consider the following family of control functions 
\begin{align*}
    h(s,t, x,y) & = u_1 |y|^{2q+16} + u_2  |x|^2 |y|^{2q + 12} + u_3 |x|^{2q + 8} \\ 
    &= (\alpha_1 \bar{y}^{q+8} + \alpha_2 \bar{x}\bar{y}^{q+6} + \alpha_3 \bar{x}^{q+4}) \partial f/\partial x  + (\alpha_1 \bar{y}^{q+8} + \beta_2 \bar{x}\bar{y}^{q+6} + \beta_3 \bar{x}^{q+4}) \partial f/\partial y\\
    &= A_1 \partial f/\partial x + A_2\partial f/\partial y
\end{align*} 
The Newton diagram of $h$ can be illustrated as in Figure \eqref{fig:newton-diagram3}

\begin{figure}[h!]
    \centering
\begin{tikzpicture}[xscale=0.4, yscale=0.4]
    \draw[->] (-1, 0) -- (15, 0) node[right] {$x$};
    \draw[->] (0, -1) -- (0, 22) node[above] {$y$};

   
 \node[below] at (4, 0) {2};
  \node[below] at (12, 0) {2q+8};

    \node[left] at (0, 20) {2q+16};
        \node[left] at (0, 12) {2q+12};

    \fill[red] (12, 0) circle (7pt) ;
    \fill[red] (4, 12) circle (7pt) ;
    \fill[red] (0, 20) circle (7pt) ;

    \draw[thick, blue] (12, 0) -- (4, 12) -- (0, 20);

    \draw[dashed, thick] (4, 12) -- (4, 0);
    \draw[dashed, thick] (4, 12) -- (0, 12);

     \node[blue] at (8, 16) {$\sigma_1 = (2, 1; 2q+16)$};

     \node[blue] at (17, 7) {$\sigma_2 = (1 + \frac{3}{q+3}, 1; 2q+14 + \frac{6}{q+3})$};

\end{tikzpicture}
\caption{Newton diagram of $h$}
\label{fig:newton-diagram3}
\end{figure}

$\Gamma(h)$ has only two compact $1$ dimensional faces $\sigma_1$ and $\sigma_2$ with weights $w_{\sigma_1} = (2, 1)$ and $w_{\sigma_2} = (1+\frac{3}{q+3}, 1)$ respectively. 
By Lemma \ref{lem_(7)} (iv), $\fil_{w_{\sigma_1}} a_i$ and  $\fil_{w_{\sigma_1} b_i}$ are $\geq \min\{\fil_{w_{\sigma_1}} x^2, \fil_{w_{\sigma_1}} xy, \fil_{w_{\sigma_1}} y^3\} = 3$. This implies that for $i = 1, 2$,
$$\fil_{w_{\sigma_1}} A_i  \geq q + 11.$$ 
In addition, $\fil_{w_{\sigma_1}} (\partial f/\partial t) = \fil_{w_{\sigma_1}} x^2 y^{4+q} = 8+q$, $\fil_{w_{\sigma_1}} h = 2q + 16$, and $w_{\sigma_1, *} = 2$.  

We have 
$$(q+11) + (8+q) - (2q + 16) -2 = 1\geq 0.$$

Thus, the condition \eqref{equ_4.1.3} in Theorem \ref{thm_checking_triviality_1} holds for $\sigma_1$. 

Consider the second face $\sigma_2$. By the same arguments, we have: 

 $\fil_{w_{\sigma_2}} \alpha_i$ and  $\fil_{w_{\sigma_2} \beta_i}$ both are $\geq$ $\min\{\fil_{w_{\sigma_2}} x^2, \fil_{w_{\sigma_2}} xy, \fil_{w_{\sigma_2}} y^3\} = 2 + \frac{3}{q+3}$.

$$\fil_{w_{\sigma_2}} A_i  \geq (q + 7 +  \frac{3}{q+3}) + (2+ \frac{3}{q+3}) = q+9 + \frac{6}{q+3}, i = 1, 2$$

and $\fil_{w_{\sigma_2}} (\partial f/\partial t) =\fil_{w_{\sigma_2}}  x^2 y^{4+q} =  6+q + \frac{6}{3+q}$, 
$\fil_{w_{\sigma_2}} h = 2q + 14+\frac{6}{q+3}$, $w_{\sigma_2,*} = 1 + \frac{3}{3+q}$. 
We have 
$$ (q + 9+ \frac{6}{q+3}) +  (6+q + \frac{6}{3+q}) - (2q + 14+\frac{6}{q+3}) - (1 + \frac{3}{3+q}) = \frac{3}{3+q}\geq 0 .$$
This implies that the condition \eqref{equ_4.1.3} in Theorem \ref{thm_checking_triviality_1} holds for $\sigma_2$. 

Therefore, $f$ is bi-Lipschitz trivial in $t$ with any $s\neq 0$ fixed. 
\end{proof}


\begin{proof}[Proof of Lemma  \ref{lem_(7)}]
(i) We first compute the partial derivatives:
\[
\begin{aligned}
\frac{\partial f}{\partial x} &= 4x^3 + 4x y^3 + 2s xy^{3+q} + 2txy^{4+q},\\[4pt]
\frac{\partial f}{\partial y} &= 6x^2 y^2 + 6y^5 + s(3+q)x^2y^{2+q} + t(4+q)x^2y^{3+q}.
\end{aligned}
\]

Computation gives
\begin{align*}
    & \left(\frac{1}{4}xy^2 +\frac{s(3+q)}{24} xy^{q+2} + \frac{t(4+q)}{24} xy^{q+3}\right)\frac{\partial f}{\partial x} - \frac{1}{6}\left(x^2 + \frac{s(q+6)}{6} y^{q+3} + \frac{t(7+q)}{6}y^{4+q}\right)\frac{\partial f}{\partial y}\\
= & 
- \frac{1}{36} y^{q+8}\left( x^{2} y^{q-3}\bigl(s(q+3)+t(q+4)y\bigr)^2
+ 6\bigl(s(q+6)+t(7+q)y\bigr)\right).
\end{align*}
This shows that $y^{8+q}$ is in $\jac(f)$.

(ii) $xy^{6+q}\in \jac(f)$.

We have 
\[
\begin{aligned}
\frac{1}{s}y^{3}\frac{\partial f}{\partial x}&-\frac{2}{3s}xy\frac{\partial f}{\partial y}+\frac{q+3}{6}y^{q+3}\frac{\partial f}{\partial x}+\frac{q+4}{6s}\,t\,y^{q+4}\frac{\partial f}{\partial x}\\[6pt]
&= \frac{2(q+6)}{3}\,x\,y^{q+6}
+ \frac{2(q+7)}{3s}\,t\,x\,y^{q+7}
+ \frac{q+3}{3}\,s\,x\,y^{2q+6}\\[4pt]
&\qquad{}+ \frac{2q+7}{3}\,t\,x\,y^{2q+7}
+ \frac{q+4}{3s}\,t^{2}\,x\,y^{2q+8}\\
& = \left( \frac{2(q+6)}{3} + \hot\right) \,x\,y^{q+6}
\end{aligned}
\]
Then, $xy^{6+q}\in \jac(f)$. This proves (ii) .

(iii) $x^{4+q}\in \jac(f)$

\medskip
\noindent\emph{Claim 1.}  $x^{q-2}y^9\in \jac(f)$ for $q\geq 4$.

\smallskip
Consider a monomial $x^n y^m$ with $m,n \geq 2$. We have
\[
x^n y^m = \tfrac{1}{6} x^{n-2} y^{m-2} \frac{\partial f}{\partial y} - x^{n-2} y^{m+3} - \frac{q+3}{6} s\, x^n y^{\,m+q} - \frac{q+4}{6} t\, x^n y^{\,m+q+1}.
\]
Thus, for $n,m\geq 2$ and $q\geq 3$,
\[
x^ny^m \in \langle \frac{\partial f}{\partial y}, x^{n-2} y^{m+3}\rangle.
\]
Replacing $(n,m)$ by $(n-2,m+3)$ repeatedly (as long as $n-2k\geq 0$), we obtain
\[
x^ny^m \in \langle \frac{\partial f}{\partial y}, x^{n-2k} y^{m+3k}\rangle,\quad k=1,2,\dots.
\]
Since (i) shows $y^{q+8}\in \jac(f)$, it follows that if $m+3k\geq q+8$ then $x^{n-2k}y^{m+3k}\in \jac(f)$ and hence $x^ny^m\in \jac(f)$.
\begin{itemize}
    \item If $n=2k$, the condition is
\[
m+3k \geq q+8 \quad\Leftrightarrow\quad \tfrac{2m+3n-16}{2}\geq q.
\]
\item  If $n=2k+1$, the condition is
\[
m+3k \geq q+8 \quad\Leftrightarrow\quad \tfrac{2m+3n-19}{2}\geq q.
\]
\end{itemize}

Applying the above, we get  $x^{q-2}y^9\in \jac(f)$ for $q\geq 4$. This proves Claim~1.

We now prove $x^{q+4}\in \jac(f)$.

Consider
\[
B := x^{1+q}\frac{\partial f}{\partial x} - \tfrac{2}{3}x^{q}y \frac{\partial f}{\partial y}.
\]
Then, 
\[
B = 4x^{q+4} - 4x^q y^{6} - \frac{2}{3}q s\,x^{q+2}y^{\,q+3} - \frac{2}{3}(q+1)t\,x^{q+2}y^{\,q+4}.
\]
Claim~1 implies $x^{2+q}y^{3+q}\in \jac(f)$ for $q\geq 3$. Hence, it remains to show $x^qy^6\in \jac(f)$.

We have 
\[
\frac{1}{6} x^{q-2} y^4 \frac{\partial f}{\partial y} 
= x^q y^6 + x^{q-2} y^9 + \frac{q+3}{6} s\, x^q y^{\,q+6} + \frac{q+4}{6} t\, x^q y^{\,q+7}.
\]
It follows that 

\[
x^q y^6 = 
 \frac{1}{6} x^{q-2} y^4 \frac{\partial f}{\partial y}  - x^{q-2} y^9 - \frac{q+3}{6} s\, x^q y^{\,q+6} - \frac{q+4}{6} t\, x^q y^{\,q+7}.
\]

The right-hand side belongs to $\langle \frac{\partial f}{\partial y}, x^{q-2}y^9, y^{q+8}\rangle$. Since $y^{q+8}\in \jac(f)$, it follows that $x^qy^6 \in \langle \frac{\partial f}{\partial x},\frac{\partial f}{\partial y}, x^{q-2}y^9\rangle$.

Now we show that $x^{q-2}y^9\in \jac(f)$ for all $q\geq 3$. Indeed, for $q=3$, we have  $x^{q-2}y^9=xy^{9} = xy^{q+3}\in \jac(f)$ by (ii); and for $q\geq 4$, $x^{q-2}y^9\in \jac(f)$ by Claim~1. 

Therefore $x^{q+4}\in \jac(f)$. 

From the construction, it is clear that  (iv) follows from (i), (ii) and (iii). The lemma is proved.
\end{proof}

\end{proof}

\section{Classification of Lipschitz unimodal germs}\label{section6}

We now present the classification of Lipschitz unimodal germs. The main result is as follows: 

\begin{thm}\label{thm_corank2}
Let $f$ be a corank $2$ isolated singularity with nonzero $4$-jet. Then, $f$ is Lipschitz unimodal if and only if it is smoothly equivalent to one of the germs in the table below:

\begin{center}
\footnotesize
\begin{longtable}{|c|c|c|c|c|}
\hline
Name & Normal form & Restrictions &  $\mu (f)$ &  $\smod (f)$ \\ 
\hline 
$J_{10}$ &  $x^3 + tx^2 y^2 + y^6$& $4t^3 + 27\neq 0$  & $10$ & \multirow{7}{*}{1} \\
$J_{2,i}$ & $ x^3 +x^2y^2 + ty^{6+i}$ & $i>0, t\neq 0$ &  $10+i$ & \\
 $W_{13}$ & $ x^4 + x y^4 + t y^6$ & &  $13$ & \\
$Z_{12}$ & $x^3 y + x y^4 + t x^2 y^3$ & & $12$ & \\
$Q_{12}$ & $x^3 + y^5 + y z^2 + t x y^4$ & & $12$ & \\
$X_{1,p}$ & $x^4 + x^2 y^2 + ty^{4+p}$ & $t\neq 0, \ p \geq 2$ & $9+p$ & \\
$Y^1_{r,s}$ & $x^{4+r} + t x^2 y^2 + y^{4+s}$ & $t\neq 0, \ r+s \geq 2$ & $9+r+s$ &\\
$E_{12}$ & $x^3 + y^7 + t x y^5$  & & $12$ &\\
$E_{13}$ &$x^3 + x y^5 + t y^8$ & &$13$ & \\
$E_{14}$ &$x^3 + y^8 + t x y^6$ & & $14$ &\\
\hline
$Z_{1,0}$        & $x^3y +  sx^2y^3 +  tx y^{6} + y^7$  & $3 s^3 + 27\neq 0$ & $15$ & \multirow{7}{*}{2}  \\
$Z_{1,1}$        & $x^3y +  x^2y^3  +  s y^{8} + ty^9$  & $ s \neq 0$ & $16$ &  \\
$W_{1,0}$        & $x^4 +  s x^2y^3  + t x^2y^4 + y^{6}$  & $ s^2 \neq 4$ & $15$ & \\
$W_{1,1}$        & $x^4 + x^2y^3  +  sy^7 + ty^8$  & $ s \neq 0$ & $16$ &  \\
$W_{1,2}$        & $x^4 + x^2y^3  +  sy^{8} + ty^9$  & $ s \neq 0$ & $17$ &  \\
$W^{\#}_{1,2q-1}$        & $(x^2 +y^3)^2  +  sxy^{4+q} + txy^{5+q}$  & $q>0, s \neq 0$ & $15+2q-1$ & \\
$W^{\#}_{1,2q}$        & $(x^2 +y^3)^2  +  sx^2y^{3+q}+ tx^2y^{4+q}$  & $q>0, s\neq 0$ & $15+2q$ &     \\
$Z_{17}$    & $x^3y + y^8 + sxy^6 + xy^7$  &   & $17$ & \\
$W_{17}$        & $x^4 + xy^5 + sy^7 + ty^8$  & & $17$ & \\
$W_{18}$        & $x^4 + y^7 + sx^2 y^4 + tx^2y^5$  & & $18$ & \\
\hline
\caption{List of Lipschitz unimodal corank $2$ germs of non-zero $4$-jets}
\label{table_corank2}
\end{longtable}
\end{center}
\end{thm}
\vspace{-1.5cm}
\begin{proof} Suppose $f$ is a function germ of corank 2 with an isolated singularity at the origin with non-zero $4$-jet. 
``$\Rightarrow$": Suppose $\Lmod(f) = 1$. If $\smod(f) = 1$, by Arnold's classification, $f$ is smoothly equivalent to one of the unimodal germs in Table \ref{table_corank2}. If $\smod(f) \geq 2$, then by Theorem \ref{thm_deform_J30}, $f$ is smoothly equivalent to one of the bimodal germs in \ref{thm_deform_J30}. 
.

``$\Leftarrow$'': Suppose $f$ is smoothly equivalent to one of the germs in Table \ref{table_corank2}. By Theorem \ref{thm_Lipschitz_simple}, $f$ is not Lipschitz simple, so $\Lmod(f) \geq 1$. We then have the following cases:

\emph{Case 1:} $\smod(f) = 1$. 

In this case, $\Lmod(f) \leq \smod(f) \leq 1$, so $\Lmod(f) = 1$.

\emph{Case 2:} $\smod(f) = 2$. 

We first prove that all germs in $Z_{1,0}$ have Lipschitz modality $1$. The argument for $W_{1,0}$ is similar. Indeed, suppose $f \in Z_{1,0}$. Then,
$$
\Lmod(f) \leq \smod(f) \leq 2.
$$
It suffices to  show that $\Lmod(f) < 2$. 

Assume, on the contrary, that $\Lmod(f) = 2$. Since, $f$ has Milnor number $15$ and the order and the Milnor number are upper semicontinuous, there exists a neighborhood $U$ of $j^k(f)$ in $J_0^k(n,1)$ such that all germs in $U$ of order $\leq 4$ and Milnor numbers $\leq 15$, for large enough $k$.

By definition, if $\Lmod(f) = 2$, then there exists a semialgebraic set $V \subset J_0^k(n,1)$ of dimension at least $4$ such that:
\begin{enumerate}[(i)]
    \item $j^k(f) \in \overline{V}$,
    \item any two distinct germs in $V$ are not bi-Lipschitz equivalent.
\end{enumerate}
Thus, every germ in $V$ must have Lipschitz modality at least $2$. As shown in Case 1, all germs with Milnor number $\mu \leq 14$ have Lipschitz modality at most $1$. Therefore, $U \cap V$ can only contain germs with Milnor number $15$. The only such candidates are those in the families $Z_{1,0}$ and $W_{1,0}$. However, by Theorem~\ref{thm_checking_triviality_main}, the Lipschitz types of germs in these families form $1$-parameter families, which contradicts condition (ii).

We have shown that germs in Table~\ref{table_corank2} with Milnor number $\leq 15$ are Lipschitz unimodal. The remaining cases can be treated similarly by induction.

Let $f$ be a germ in Table~\ref{table_corank2} with Milnor number $\mu(f) > 15$. By Lemma~\ref{lem_main_deform}, $f$ does not deform to $J_{3,0}$. Thus, for $k$ sufficiently large, there exists a neighborhood $U \subset J_0^k(n,1)$ of $j^k(f)$ containing no germs that deform to $J_{3,0}$. Moreover, by upper semicontinuity of the Milnor number, any germ in $U$ has Milnor number at most $\mu(f)$.

Suppose, on the contrary, that $\Lmod(f) = 2$. Then, as before, there exists a semialgebraic set $V \subset J_0^k(n,1)$ of dimension at least $4$ such that the conditions (i) and (ii) above hold. As in the previous case, the only germs possibly contained in $V$ are those of Milnor number $\mu(f)$ from Table~\ref{table_corank2}. By Theorem~\ref{thm_checking_triviality_main}, these belong to finitely many families whose Lipschitz types are given by $1$-parameter families. Hence, condition (ii) fails. Therefore, $\Lmod(f) < 2$, completing the proof.
\end{proof}

A direct consequence of the above theorem is the following:

\begin{cor}
The germs in the family $J_{3,0}$ are Lipschitz bimodal and have the smallest possible Milnor number among all Lipschitz bimodal germs.
\end{cor}

\begin{cor}
Let $f$ be a corank $2$ germ with nonzero $4$-jet. Then, $f$ is Lipschitz unimodal if and only if it deforms to $J_{10}$ but does not deform to $J_{3,0}$.
\end{cor}

The following result provides an upper bound on the Lipschitz modality of function germs.

\begin{prop}
All isolated singularities with zero $6$-jet deform to $J_{3,0}$. Consequently, they have Lipschitz modality at least $2$.
\end{prop}

\begin{proof}
Since any zero $6$-jet of corank $\geq 2$ can deform to a zero $6$-jet of corank $2$, it suffices to consider corank $2$ singularities with vanishing $6$-jets.

Let $f(x,y)$ be an isolated singularity with $j^6(f) = 0$. Then $f$ admits an expansion of the form
$$
f(x,y) = a_0 x^7 + a_1 x^6 y + \dots + a_7 y^7 + b_1 x^8 + \dots + b_8 y^8 + \text{h.o.t.}
$$

Consider the deformation
$$
F_t(x,y) = f(x,y) + t x^7.
$$
Then
$$
F_t(x,y) = (a_0 + t) x^7 + a_1 x^6 y + \dots + a_7 y^7 + b_1 x^8 + \dots + b_8 y^8 + \text{h.o.t.}
$$

By applying a coordinate change of the form $x \mapsto x + \alpha y$, for a suitable $\alpha$, we can eliminate the $y^7$ term. That is, for $t \neq 0$ sufficiently small, we have
$$
F_t(x,y) \sim_{\al R} G_t(x,y),
$$
where $j^6(G_t) = 0$ and $G_t$ contains no $y^7$ term.

Now consider the deformation
$$
H_{t,s}(x,y) = G_t(x,y) + s x^4.
$$
For $s \neq 0$ close to $0$, we can eliminate the $y^8$ term via a coordinate change of the form $x \mapsto x + \beta y^2$, for a suitable $\beta$. It is easy to check that the resulting family has filtration $\geq 9$ with respect to weight $w = (3,1)$. By Lemma~\ref{lem_main_deform}, $H_{t,s}$ deforms to $J_{3,0}$.
Hence, $f(x,y)$ deforms to $J_{3,0}$ as claimed.
\end{proof}

\section{Final remarks and Open questions}\label{section7}
It is clear from the definition that if $ f \sim_{\mathcal{R}} g $, then $ \Lmod(f) = \Lmod(g) $.

\begin{question}
    Suppose $ f \sim_{\mathrm{Lip}} g $. Is it true that $ \Lmod(f) = \Lmod(g) $?
\end{question}
In fact, it follows directly from the  classification of Lipschitz simple germs in~\cite{nhan1} that if $f \sim_{Lip} g $ and $\Lmod(f) = 0$, then $\Lmod(g) = 0$ as well.

The next question is motivated by Theorem~\ref{thm_Lipschitz_simple2}.

\begin{question}
    Let $ f \in \mathfrak{m}_n^2 $ be a germ with an isolated singularity at the origin. Is it true that $ f $ is Lipschitz unimodal if and only if it deforms to $ J_{10} $ but does not deform to $ J_{3,0} $?
\end{question}

It is well known that Thom's slitting lemma plays a fundamental role in the classification theory of singularities. However, it remains unclear whether a Lipschitz version of this lemma holds. More precisely:

\begin{question}
    Let $ f, g \in \mathfrak{m}_n^2 $ be germs with isolated singularities at the origin. Suppose that 
    $$
    f(x) + Q(z) \sim_{\Lip} g(x) + Q(z),
    $$
    where $ Q(z) = z_1^2 + \dots + z_m^2 $. Does it follow that $ f \sim_{\mathrm{Lip}} g $?
\end{question}

In \cite[Theorem 5.1]{nhan1} It was shown that under the above assumption, $f$ and $g$ must have the same multiplicity. Moreover, their principal homogeneous parts are bi-Lipschitz equivalent. 

\section*{Acknowledgement}
This work was initiated in collaboration with Prof. Maria Aparecida Soares Ruas. We note with deep sadness that she passed away before its completion and gratefully acknowledge her invaluable ideas, guidance and inspiration. The first author would like to express his gratitude to the Vietnam Institute for Advanced Study in Mathematics (VIASM) for its warm hospitality and generous support during the preparation of this paper. 
We thank H.D.~Nguyen for his interest and valuable discussions. 
The first author was supported by the Vietnam Ministry of Education and Training (MOET) under grant number B2025-CTT-01.

\end{document}